\documentclass[12pt]{article}
\makeatletter
\renewcommand{\@makefnmark}{}
\makeatother

\usepackage{amscd}
\usepackage{amsmath,amssymb}
\usepackage{amssymb}
\usepackage{amsmath}
\usepackage{amsthm}
\usepackage{latexsym}
\usepackage{amsfonts}
\usepackage{color}
\usepackage{graphics}
\usepackage[all,knot]{xy}
\usepackage[all]{xy}
\xyoption{arc}
\usepackage{enumerate}

\usepackage{lipsum}

\usepackage{a4wide}
\usepackage{eurosym}
\usepackage{makeidx}

\usepackage{amstext}
\usepackage{secdot}

\theoremstyle{plain}
\newtheorem{thm}{\indent\sc Theorem}[section]
\newtheorem{prop}[thm]{\indent\sc Proposition}

\newtheorem{cor}[thm]{\indent\sc Corollary}
\theoremstyle{definition}
\newtheorem{defn}[thm]{\indent\sc Definition}

\newtheorem{example}[thm]{Example}
\theoremstyle{remark}

\def\proofwp{{\sc Proof}}

\def\ex{\begin{example}\rm}
\def\endex{\qed\end{example}}
\def\rem{\begin{rem}\rm}
\def\endrem{\qed\end{rem}}

\def\qed{\hfill{$\square$}}

\def\mapr#1{\smash{\mathop{\buildrel{#1}\over\longrightarrow}}}

\def\mapd#1{\Big\downarrow\rlap{$\vcenter{\hbox{$#1$}}$}}

\def\mapr#1{\smash{\mathop{\buildrel{#1}\over\longrightarrow}}}

\def\mapd#1{\Big\downarrow\rlap{$\vcenter{\hbox{$#1$}}$}}

\def\qed{\hfill\vrule width2mm height2mm depth2mm}
\def\v|{\;\vrule width1pt height3mm depth0mm\;}

\def\H{{\bf H}}

\def\Z{{\bf Z}}

\def\s{{\bf s}}

\def\0{{\bf 0}}
\def\1{{\bf 1}}

\def\cA{{\cal A}}
\def\cB{{\cal B}}
\def\cC{{\cal C}}
\def\cD{{\cal D}}

\def\cO{{\cal O}}

\def\ad{{\hbox{\bf ad}}}

\def\Out{{\hbox{\bf Out}\;}}

\def\End{{\hbox{\bf End}}}

\def\id{{\hbox{\bf Id}}}
\def\im{{\hbox{\bf Im}\;}}

\def\ker{{\hbox{\bf Ker}\;}}

\def\supp{{\hbox{\bf supp}\;}}

\def\<{\langle}
\def\>{\rangle}
\def\({\left(}
\def\){\right)}

\def\St{\mbox{\bf St }}

\newcommand{\RR}{\mathbb{R}}

\def\beq#1{\color{blue}\begin{equation}\label{#1}}
\def\eeq{\end{equation}\color{black}}

\title{{\bf Sullivan constructions for transitive Lie algebroids - smooth case}}
\author{{\bf Aleksandr S. Mishchenko} \\
{\bf Jose R. Oliveira}}

\newcommand{\Addresses}{{
  \bigskip
  \footnotesize

  A.S.~Mishchenko, \textsc{Department of high Geometry and Topology, Moscow State University,
    Moscow, Russia}\par\nopagebreak
  \textit{E-mail address}: \texttt{asmish@mech.math.msu.su}

  \medskip

  J.R.~Oliveira, \textsc{Department of Mathematics and Applications, University of Minho,
    Braga, Portugal}\par\nopagebreak
  \textit{E-mail address}: \texttt{jmo@math.uminho.pt}

}}

\date{}

\begin{document}

\maketitle

\pagenumbering{arabic}

\begin{center}

\vspace{3mm}

{\large \textbf{Abstract}}

\vspace{3mm}

\end{center}

\addcontentsline{toc}{section}{\hspace*{0.5cm}Abstract}


Let $M$ be a smooth manifold, smoothly triangulated by a simplicial complex $K$, and $\cA$ a transitive Lie algebroid on $M$. The Lie algebroid restriction of $\cA$ to a simplex $\Delta$ of $K$ is denoted by $\cA^{!!}_{\Delta}$. A piecewise smooth form of degree $p$ on $\cA$ is a family $\omega=(\omega_{\Delta})_{\Delta\in K}$ such that $\omega_{\Delta}\in \Omega^{p}(\cA^{!!}_{\Delta};\Delta)$ for each $\Delta\in K$, satisfying the compatibility condition concerning the restrictions of $\omega_{\Delta}$ to the faces of $\Delta$, that is, if $\Delta'$ is a face of $\Delta$, the restriction of the form $\omega_{\Delta}$ to the simplex $\Delta'$ coincides with the form $\omega_{\Delta'}$. The set $\Omega^{\ast}(\cA;K)$ of all piecewise smooth forms on $\cA$ is a cochain algebra. One has a natural morphism
$$\Omega^{\ast}(\cA;M)\rightarrow \Omega^{\ast}(\cA;K)$$
of cochain algebras given by restriction of a smooth form defined on $\cA$ to a smooth form defined on $\cA^{!!}_{\Delta}$, for all simplices $\Delta$ of $K$. In this paper, we prove that, for triangulated compact manifolds, the cohomology of this construction is isomorphic to the Lie algebroid cohomology of $\cA$, in which the isomorphism is induced by the restriction map.

\begin{center}
\tableofcontents
\end{center}

\vspace{3mm}

\begin{center}

\vspace{3mm}
{\large \textbf{Introduction}}
\end{center}

\addcontentsline{toc}{section}{\hspace*{0.5cm}Introduction}

\vspace{3mm}

D. Sullivan considered in \cite{suli-inf} a new model for the underlying cochain complex of classical cohomologies with rational coefficients for arbitrary simplicial spaces which gives an isomorphism with classical rational cohomologies. This new model is determined by the Rham complex of all rational polynomial forms defined on the simplicial complex triangulating the space. H. Whitney also presented in \cite{wity-git} other cell-like constructions of cochain complexes which induce isomorphisms in cohomology with classical cohomologies.

The present paper arose from efforts to extend those constructions to transitive Lie algebroids and to apply Sullivan's methods for studying
formality of transitive Lie algebroids. Among the constructions presented in \cite{suli-inf} and \cite{wity-git}, there is one which states that de Rham cohomology of a smooth manifold smoothly triangulated by a simplicial complex is isomorphic to piecewise smooth cohomology of the simplicial complex. This isomorphism is given by restriction of smooth forms to all simplices.

The key ideas concerning class obstruction arising from non-abelian extensions of Lie algebroids have inspired us and led to conjecture that,
given a transitive Lie algebroid on a triangulated compact smooth manifold, the morphism given by restriction, which
takes smooth forms on the Lie algebroid into piecewise smooth forms on the same Lie algebroid, still remains an
isomorphism in cohomology.

The aim of the present paper is to prove this conjecture. For this purpose, we use the structure which commences by fixing a smooth triangulation of the base of a transitive Lie algebroid by a simplicial complex and taking the restriction of the Lie algebroid to all simplices of the triangulation. Since the Lie algebroid is transitive, the restriction of the Lie algebroid to each simplex always exists. When this structure is given, we define the notion of piecewise smooth form in a similar way to Whitney forms on a simplicial complex. The set of all piecewise smooth forms defined on a transitive Lie algebroid over a triangulated base is naturally equipped with a differential, yielding a commutative differential graded algebra. Its cohomology is, by definition, the piecewise smooth cohomology of the Lie algebroid. Each smooth form defined on the Lie algebroid gives a piecewise smooth form defined by taking the restriction of the form to each simplex. This correspondence is a natural map from the usual algebra of the smooth forms on the Lie algebroid to the algebra of the piecewise smooth forms of the same Lie algebroid. Based on three crucial results, namely the triviality of a transitive Lie algebroid over a contractible smooth manifold (Mackenzie), the K$\ddot{\textrm{u}}$nneth theorem for Lie algebroids (Kubarski) and the de Rham-Sullivan theorem for smooth manifolds, we show that this map is an isomorphism in cohomology.

We have tried to make this paper as self-contained as possible. It is for this reason that some parts of this paper include basic constructions of Lie algebroids and a treatment on non-abelian extensions of Lie algebroids. Throughout the paper, all manifolds are smooth, finite-dimensional and possibly with boundaries of different indices.

\vspace{3mm}

\textbf{Acknowledgments}. We want to thank to James Stasheff for his strong dynamism to discuss several topics concerning this work. We also thank to Nicolae Teleman for his valuable ideas on theory of Lie algebra extensions. This paper is dedicated to Jan Kubarski who already departed but was always with us, giving his wholehearted support in Chern-Weil theory for Lie algebroids. The second author much appreciates the support given by Paula Smith, Lisa Santos and Fernando Miranda.

\vspace{6mm}

\begin{center}
\section{Restrictions of transitive Lie algebroids}
\end{center}

\vspace{3mm}

The restriction of a Lie algebroid to an open subset of the base coincide with the restriction of the underlying vector bundle. For more general submanifolds of the base, the restricted Lie algebroid may not coincide with the restriction of the underlying vector bundle, but it will be a vector subbundle of the
restriction of the underlying vector bundle. Restrictions of Lie algebroids to submanifolds of the base are given by inverse image of Lie algebroids through the inclusion of submanifolds. The transitivity of Lie algebroids plays an important role on restrictions to general submanifolds. In order to have a coherent notation and terminology throughout the paper, we shall summarize some definitions and examples of Lie algebroids. A reasonably complete description of definitions and results on Lie algebroids can be found in Mackenzie \cite{makz-lga}, \cite{makz-gad} and \cite{makz-acla} and Kubarski \cite{kuki-chw}. A short introduction to the categorical point of view of Lie algebroids can be seen in the preprint $``$\textit{Transitive Lie algebroids - categorical point of view}$"$ by Mishchenko \cite{mish-tla}.

\vspace{3mm}

Let $M$ be a smooth manifold, possibly with boundary and corners, $TM$ the tangent bundle to $M$ and $\Gamma(TM)$ the Lie algebra of the vector
fields on $M$. We recall that a Lie algebroid on $M$ is a vector bundle $\pi:\cA\longrightarrow M$ on $M$ equipped with a vector bundle
morphism $\gamma:\cA\longrightarrow TM$, called anchor of $\cA$, and a structure of real Lie algebra on the vector space $\Gamma(\cA)$
of the sections of $\cA$ such that the map $\gamma_{\Gamma}:\Gamma(\cA)\mapr{}\Gamma(TM)$, induced by $\gamma$, is a Lie algebra homomorphism
and the action of the algebra $\cC^{\infty}(M)$ on $\Gamma(\cA)$ satisfies the natural condition:
$$[\xi,f\eta]=f[\xi,\eta] + (\gamma_{\Gamma}(\xi)\cdot f)\eta$$ for each $\xi$, $\eta$ $\in \Gamma(\cA)$ and
$f\in \cC^{\infty}(M)$. The Lie algebroid $\cA$ is called transitive if the anchor $\gamma$ is surjective. As usually, when there is no ambiguity,
we drop the anchor map and the Lie bracket in the notation of the Lie algebroid but, when it is needed to emphasize them,
we write $(\cA,[\cdot,\cdot],\gamma)$ for denoting this structure.

If $(\cB,[\cdot,\cdot],\delta)$ is other Lie algebroid on a smooth manifold $N$, a morphism of Lie algebroids from $(\cA,[\cdot,\cdot],\gamma)$
to $(\cB,[\cdot,\cdot],\delta)$ consists of a pair of mappings $(\psi,\varphi)$, with $\psi:\cA\longrightarrow \cB$ and $\varphi:M\longrightarrow N$,
such that $(\psi,\varphi)$ is a vector bundle morphism satisfying the equality $\delta\circ \psi=T(\varphi)\circ \gamma$, where
$T(\varphi):TM\longrightarrow TN$ means the tangential of $\varphi$, and preserving the Lie bracket condition for $\psi$-decompositions,
that is, for each $\xi,\eta\in \Gamma(\cA)$ with decompositions
$$\psi\circ \xi=\sum_{i=1}^{m}a_{i}\otimes \xi_{i} \ \ \ \ \ \ \ \ \ \psi\circ \eta=\sum_{j=1}^{m}b_{j}\otimes \eta_{j}$$
then
$$\psi\circ [\xi,\eta]=\sum_{i,j}a_{i}b_{j}\otimes [\xi_{i},\eta_{j}]+\sum_{j=1}^{m}(\gamma\circ\xi)(b_{j})\otimes\eta_{j}-\sum_{i=1}^{m}(\gamma\circ\eta)(a_{i})\otimes\xi_{j}$$
When $M=N$ and $\varphi$ is the identity map, the Lie algebroid morphism is called a strong Lie algebroid morphism and denoted simply by $\psi$.
We notice that, in this case, a simple characterization to ensure that a vector bundle morphism between two Lie algebroids is a Lie algebroid morphism
can be seen in \cite{kuki-chw} or \cite{makz-lga}. Namely, if $(\cA,[\cdot,\cdot],\gamma)$ and $(\cB,[\cdot,\cdot],\delta)$ are Lie algebroids
over the same smooth manifold $M$, then a vector bundle morphism $\psi$ from $\cA$ to $\cB$ is a Lie algebroid morphism if, and only if,
$\gamma=\psi\circ \delta$ and the map $\psi_{\Gamma}:\Gamma(\cA)\longrightarrow \Gamma(\cB)$, induced by $\psi$, is a Lie algebra morphism.

\vspace{3mm}

We give some examples of Lie algebroids. The first three examples are widely used in this work.

\vspace{3mm}

\textit{Example 1} (Lie algebras). Any real finite dimensional Lie algebra $\frak g$ over a one-point space $M=\{\ast\}$
with anchor equal to zero is a Lie algebroid on $M$. Any Lie algebra morphism between two Lie algebras is a Lie algebroid morphism
for this structure of Lie algebroid.

\vspace{3mm}

\textit{Example 2} (Tangent Lie algebroid). If $M$ is a smooth manifold then $TM$ is a Lie algebroid on $M$. The anchor map
is the identity map of $TM$, and the Lie bracket is the usual Lie bracket of vector fields. This Lie algebroid is called the tangent Lie algebroid of $M$. More generally, if $F$ is a regular foliation in $M$, then the tangent Lie algebroid of $F$ is the vector subbundle of $TM$ consisting
of tangent spaces to $F$, with the usual Lie bracket, and the inclusion map as the anchor. Any Lie algebroid whose anchor map is injective is isomorphic to the tangent algebroid of some regular foliation (see Kubarski \cite{kuki-chw}). Let $(\cA,[\cdot,\cdot],\gamma)$ be a Lie algebroid over a smooth manifold M. Then, the anchor map $\gamma:\cA\longrightarrow TM$ is a Lie algebroid morphism from $\cA$ to the tangent Lie algebroid of $M$.

\vspace{3mm}

\textit{Example 3} (Trivial Lie algebroid). Let $\frak g$ be a real finite dimensional Lie algebra and $M$ a smooth manifold. Consider the trivial vector bundle $M\times \frak g$ of base $M$. On the fibre product $$TM\oplus (M\times \frak g)=TM\times_{M} (M\times\frak g)$$ we define an anchor map $\gamma:TM\oplus (M\times \frak g)\longrightarrow TM$ by taking $\gamma$ to be the projection of $TM\oplus (M\times \frak g)$ on $TM$ and a Lie bracket on $\Gamma (TM\oplus (M\times \frak g))$ by setting $$[\big(X,u),(Y,v)]=([X,Y]_{TM}, X(v)-Y(u)-[u,v]\big)$$ for $X,Y\in \Gamma (TM)$ and $u,v:M\longrightarrow \frak g$ smooth maps. The vector bundle $TM\oplus (M\times \frak g)$, equipped with this anchor and this Lie bracket on the set of the sections, is a Lie algebroid on $M$ and is called trivial Lie algebroid on $M$ with fibre $\frak g$.

\vspace{3mm}

\textit{Example 4} (Lie algebroid product). Let $M$ and $N$ be two smooth manifolds and $(\cA,[\cdot,\cdot]_{\cA},\gamma)$ and $(\cB,[\cdot,\cdot]_{\cB},\widehat{\gamma})$ Lie algebroids over $M$ and $N$ respectively. The product of the Lie algebroids $(\cA,[\cdot,\cdot]_{\cA},\gamma)$ and $(\cB,[\cdot,\cdot]_{\cB},\widehat{\gamma})$, denoted by $\cA\times \cB$, is the vector bundle product $\cA\times \cB$ over $M\times N$, in which the anchor is $\gamma\times \widehat{\gamma}$ and the Lie bracket is defined in the following way: for each $\xi=(\xi_{1},\xi_{2})$ and $\eta=(\eta_{1},\eta_{2})$ $\in \Gamma(\cA\times \cB)$ $$[\xi,\eta]_{\cA\times \cB}=([\xi,\eta]^{1},[\xi,\eta]^{2})\in \Gamma(\cA\times \cB)$$ where
$$[\xi,\eta]^{1}_{(x,y)}=[\xi_{1}(-,y),\eta_{1}(-,y]_{\cA}(x)+\widehat{\gamma}(\xi_{2_{(x,y)}})(\eta_{1}(x,-))
-\widehat{\gamma}(\eta_{2_{(x,y)}})(\xi_{1}(x.-))$$ and
$$[\xi,\eta]^{2}_{(x,y)}=[\xi_{2}(-,y),\eta_{2}(-,y]_{\cA}(x)+\gamma(\xi_{1_{(x,y)}})(\eta_{2}(x,-))
-\gamma(\eta_{1_{(x,y)}})(\xi_{2}(x.-))$$
Consider now a real finite dimensional Lie algebra $\frak g$ and the tangent Lie algebroid $TM$. The Lie algebra $\frak g$ is a Lie algebroid over a one-point space $N=\{\ast\}$. We may take the product of Lie algebroids $TM \times \frak g$. The Lie algebroid $TM \times \frak g$ is defined over $M\simeq M\times N$. On the other hand, we can consider the trivial Lie algebroid $TM\oplus (M\times \frak g)$ over $M$ and we easily see that the map $f:TM\times \frak g\longrightarrow TM\oplus (M\times \frak g)$ given by $f(x,u,v)=(x,u,x,v)$ is an isomorphism of Lie algebroids. Henceforth, we will identify both Lie algebroids.

\vspace{3mm}

\textit{Example 5} (Lie algebra bundle). A Lie algebra bundle over a smooth manifold $M$ is a vector bundle
$\pi:L\longrightarrow M$ equipped with a section $[\cdot,\cdot]$ of the vector bundle $\bigwedge^{2}(L,L)$ such that,
for each $x\in M$, $(L_{x},[\cdot,\cdot]_{x})$ is a real Lie algebra and $L$ admits an atlas
$$\{\psi_{j}:U_{j}\times \frak g\longrightarrow \pi^{-1}(U), j\in J\} \ \ \ \ \ (\frak g \ \textrm{is a Lie algebra})$$
in which each $\psi_{j_{x}}$ is a Lie algebra isomorphism. Any Lie algebra bundle is a Lie algebroid, in which its anchor is equal to zero
and the Lie bracket is defined in natural way. We show now two examples of Lie algebra bundles.

Let $\pi:E\longrightarrow M$ be a vector bundle on a smooth manifold $M$. Then, the vector bundle $\End(E)=L(E;E)$, whose fibres,
at each point $x\in M$, are the vector spaces $L(E_{x};E_{x})$, is a Lie algebra bundle (see \cite{makz-lga}).

Let $\cA$ be a Lie algebroid on $M$, with anchor $\gamma:\cA\longrightarrow TM$, and consider the Atiyah sequence
$$
\xymatrix{
\{0\} \ar[r] & \ker \gamma \ar[r]^{j} & \cA \ar[r]^{\gamma} & TM \ar[r] & \{0\}\\
}
$$ Then, $\ker \gamma$ is a Lie algebra bundle on $M$.

\vspace{3mm}

\textit{Example 6} (Lie algebroid of a Lie groupoid). Let $G$ and $M$ two smooth manifolds. The manifold $G$ is called a Lie groupoid on $M$
if it is given two surjective submersion $\alpha:G\longrightarrow M$ and $\beta:G\longrightarrow M$, called the source projection and
the target projection respectively, a smooth map $1:M\longrightarrow G$ called the object inclusion map, a smooth multiplication in
$G\ast G=\{(h,g)\in G\times G: \alpha(h)=\beta(g)\}$ and a map $G\longrightarrow G$ called inverse map and denoted by $g\longrightarrow g^{-1}$,
satisfying certain identities (for the details, see Mackenzie \cite{makz-lga}, definition 1.1.1 and definition 1.1.3). Denote by $\cA(G)$ the vector bundle $\bigsqcup _{x\in M}T_{1_{x}}G_{x}$ (disjoint union). Consider the map $\gamma:\cA(G)\longrightarrow TM$ defined by $\gamma(a)=D\beta_{1_{x}}(a)$. Define a Lie bracket on $\Gamma(\cA(G))$ in the following way: for each $\xi$ and $\eta$ $\in \Gamma(\cA(G))$, the Lie bracket is defined by
$$[\xi,\eta]=[\xi',\eta']_{G}$$ in which $\xi'$ and $\eta'$ denote the unique $\alpha$-right-invariant vector fields on $G$ such that
$\xi'_{1_{x}}=\xi_{x}$ and $\eta'_{1_{x}}=\eta_{x}$, $\forall x\in M$. Then, $(\cA(G),[\cdot,\cdot],\gamma)$ is a Lie algebroid on $M$
and is called the Lie algebroid of the Lie groupoid $G$. Lie groupoids and Lie algebroids enjoy some of the properties of Lie groups and Lie algebras.
We notice that not every Lie algebroid is integrable to a Lie groupoid. The theorem 4.1 of the paper \cite{crac-loja} shows necessary and sufficient conditions so that a Lie algebroid is integrable to a Lie groupoid.

\vspace{3mm}

An example of a Lie algebroid constructed from a Lie groupoid is the Lie algebroid of the covariant derivatives of a Lie algebra bundle.
We describe this example now. Let $M$ be a smooth manifold and $\pi:E\longrightarrow M$ a vector bundle on $M$. Denote by $\Phi(E)$ the Lie groupoid on $M$ made by all linear isomorphism $\xi:E_{x}\longrightarrow E_{y}$ for each $x$, $y$ $\in M$, in which the source map $\alpha:\Phi(E)\longrightarrow M$ and the target projection $\beta:\Phi(E)\longrightarrow M$ are defined by $\alpha(\xi)=x$ and $\beta(\xi)=y$, the object inclusion map
$1:M\longrightarrow \Phi(E)$ is defined by $x\longrightarrow id_{E_{x}}$ and the multiplication is the composition of maps (see example 1.1.12 of \cite{makz-lga}). For each $n\geq 1$, let $\pi_{L}:L^{n}(E;E)\longrightarrow M$ be the vector bundle whose fibre at $z\in M$ is the vector space of all $p$-linear maps from $E_{z}\times\dots\times E_{z}$ to $E_{z}$. The canonical action $$\Phi(E)\ast L^{n}(E;E)\longrightarrow L^{n}(E;E)$$ of the Lie groupoid $\Phi(E)$ on the vector bundle $\pi_{L}:L^{n}(E;E)\longrightarrow M$ is defined by
$$\xi\cdot \varphi=\xi\circ \varphi\circ (\xi^{-1}\times\dots\times \xi^{-1})\in L^{n}(E_{y};E_{y})$$
where $x$, $y$ $\in M$, $\xi:E_{x}\longrightarrow E_{y}$ is a linear isomorphism and $\varphi\in L^{n}(E_{x};E_{x})$
(see section 1.6 of \cite{makz-lga} for more details on actions of Lie groupoids). A section $\eta \in \Gamma (L^{n}(E;E))$ is stable for this action if, for all $x$, $y$ $\in M$, there is a linear isomorphism $\xi:E_{x}\longrightarrow E_{y}$ such that $\xi\cdot \eta_{x}=\eta_{y}$.
For a stable section $\eta \in \Gamma (L^{n}(E;E))$, the stabilizer of $\Phi(E)$ at $\eta$ is defined by
$$\{\xi\in \Phi(E): \ \ \xi\cdot \eta(\alpha(\xi))=\eta(\beta(\xi))\}$$ The stabilizer of $\Phi(E)$ at $\eta$ is a Lie groupoid on $M$
(see example 1.3.4 and theorem 1.6.23 of \cite{makz-lga}). We notice that a section $\eta \in \Gamma (L^{n}(E;E))$ need not be stable. Nevertheless, for a Lie algebra bundle $\pi':K\longrightarrow M$ on $M$ with bracket $[\cdot,\cdot]\in \bigwedge^{2}(K,K)$, the bracket $[\cdot,\cdot]$ is a stable section for the action above restricted to the vector bundle $\bigwedge^{2}(K;K)$ (example 1.7.12 of \cite{makz-lga}). Hence, the stabilizer of $\Phi(K)$ at $[\cdot,\cdot]$ is well defined (theorem 1.6.23 of \cite{makz-lga}). Its Lie algebroid is denoted by $\cD_{Der}(K)$ and called the \textit{Lie algebroid of the covariant derivatives} of $K$. The Lie algebroid $\cD_{Der}(K)$ is transitive (see \cite{makz-lga}, \cite{mish-cla}, \cite{mish-latb}, \cite{mish-obs}).

\vspace{3mm}

\textit{Example 7} (Adjoint Lie algebra bundle). Let $M$ be a smooth manifold and $\pi:K\longrightarrow M$ a Lie algebra bundle on $M$ with fibre type $\frak g$. Consider the Lie subalgebra $\textbf{Der}(\frak g)$ of the Lie algebra $\frak gl(\frak g)$ made by the derivations of $\frak g$. By the proposition 3.3.9 from \cite{makz-lga}, it can be seen that the Lie subalgebra $\textbf{Der}(\frak g)$ corresponds to a unique Lie algebra subbundle of the Lie algebra bundle $\End(K)$. This Lie algebra subbundle is denoted by $\textbf{Der}(K)$ and its elements are called \textit{derivations} of $K$. Thus, the image $\ad(K)$ of $K$, by the Lie algebra bundle morphism $\ad:K\longrightarrow \textbf{Der}(K)$ defined by $\ad_{x}:K_{x}\longrightarrow \textbf{Der}(K_{x})$ for each $x\in M$, is a Lie algebra subbundle of $\textbf{Der}(K)$, in which the fibre $\ad(K)_{x}$ is an ideal of $\textbf{Der}(K_{x})$. Consequently, the Lie algebra bundle quotient $\textbf{Der}(K)/\ad(K)$ is well defined (see \cite{makz-lga}). The Lie algebra subbundle $\ad(K)$ is called the adjoint Lie algebra bundle of $K$. The Lie algebra bundle quotient $\textbf{Der}(K)/\ad(K)$ usually is denoted by $\Out(K)$.

\vspace{3mm}

Assume now that $\cA$ is a transitive Lie algebroid over $M$, with anchor $\gamma:\cA\longrightarrow TM$ and Lie bracket $[\cdot,\cdot]_{\cA}$ on $\Gamma (\cA)$.
Consider the Lie algebra bundle $\ker \gamma$. We recall that an ideal of $\cA$ is a Lie algebra subbundle $K$ of $\ker \gamma$ such that,
for all sections $\xi\in \Gamma (\cA)$ and $\eta\in \Gamma K$, $[\xi,\eta]_{\cA}$ is a section of $K$. In these conditions, one defines the
\textit{Lie algebroid quotient} of $\cA$ by $K$ as follows. Let $\overline{\cA}$ be the vector bundle quotient $\cA/K$ and
$\overline{\gamma}:\overline{\cA}\longrightarrow TM$ the map induced by $\gamma$. The Lie bracket in the space of the sections of
$\overline{\cA}$ is defined by $$[\xi+\Gamma K,\eta+\Gamma K]_{\overline{\cA}}=[\xi,\eta]_{\cA}+\Gamma K$$ for each $\xi$, $\eta$ $\in \Gamma (\cA)$.
The Lie algebroid $\overline{\cA}$ is transitive and usually denoted by $\cA/K$ (see proposition 6.5.8 of \cite{makz-lga} or proposition IV - 1.11 of \cite{makz-gad} for more details). A particular case of Lie algebroid quotient is the following example. Let $\pi:K\longrightarrow M$ be a Lie algebra bundle on a smooth manifold $M$ and consider the transitive Lie algebroid $\cD_{Der}(K)$ of all covariant derivatives of a Lie algebra bundle $K$. The adjoint Lie algebra bundle $\ad(K)$ is an ideal of $\cD_{Der}(K)$. Hence, we can consider the transitive Lie algebroid quotient $\cD_{Der}(K)/\ad(K)$, which is denoted by $\Out\cD(K)$. The sections of the Lie algebroid $\Out\cD(K)$ are called outer derivations of $K$.

\vspace{3mm}

Let $M$ be a smooth manifold and $\cA$ a Lie algebroid over $M$, with anchor $\gamma:\cA\longrightarrow TM$ and Lie bracket $[\cdot,\cdot]$ on $\Gamma (\cA)$.
We recall that the Lie bracket $[\cdot,\cdot]$ enjoys the local property, that is, if $U$ is an open subset of $M$ and $\xi,\eta\in \Gamma (\cA)$ such that $\eta$ vanishes on $U$ then the Lie bracket $[\xi,\eta]$ vanishes on $U$. This makes the restricted vector bundle $\cA_{U}$ into a Lie algebroid on $U$ (see Mackenzie \cite{makz-lga}, section 3.3). Under suitable conditions, we can find a rather understandable and complete notion of restrictions of the Lie algebroid $\cA$ to more general submanifolds of $M$. Our definition of restriction is based on the construction of inverse image of a Lie algebroid by a smooth map. The transitivity of the Lie algebroid $\cA$ will be needed for showing that inverse image always exists for any smooth map. We notice now some brief considerations on the construction of image inverse.

Assume that $\cA$ is transitive. Let $N$ be other smooth manifold and $\varphi:N\longrightarrow M$ a smooth map. Let $\pi:\cA\longrightarrow M$ denote the vector bundle underlying the Lie algebroid $\cA$ and $\pi_{TM}:TM\longrightarrow M$ and $\pi_{TN}:TN\longrightarrow N$ the canonical projections. The anchor $\gamma:\cA\longrightarrow TM$ defines $\cA$ as a vector bundle on $TM$. In order to complete the following diagram
\begin{equation*}
\begin{array}{ccc}
\  &  & \cA \\
& & \mapd{\gamma}\\
\\ TN & \mapr{T\varphi} & TM
\end{array}
\end{equation*}
consider the vector bundle $\widehat{\gamma}:(T\varphi)^{\ast}\cA\longrightarrow TN$, which is the inverse image of the vector bundle $\gamma:\cA\longrightarrow TM$ by the smooth map $T\varphi$, where $\widehat{\gamma}$ denotes the canonical projection $\widehat{\gamma}:(T\varphi)^{\ast}\cA\longrightarrow TN$. We notice that the vector bundle $\widehat{\gamma}:(T\varphi)^{\ast}\cA\longrightarrow TN$ exists because $\gamma$ is surjective. Obviously, $(T\varphi)^{\ast}\cA$ is also a vector bundle over $N$, in which its projection is the composition of the projection $\widehat{\gamma}$ with the canonical projection $\pi_{TN}$. We obtain the following commutative diagram
$$
    \xymatrix@R20pt@C20pt{&
      (T\varphi)^{*}(\cA)\ar[r]\ar[d]_{\widehat{\gamma}} &  \cA\ar[d]_{\gamma} \ar@/^2pc/[dd]^{\pi}\\
      &TN\ar[r]^{T\varphi}\ar[d]_{\pi_{N}}&TM\ar[d]_{\pi_{M}}\\
      &N \ar[r]^{\varphi}&M
    }
$$
Consider the Whitney sum $TN\oplus\varphi^{\ast}(\cA)$ of the tangent bundle $TN$ on $N$ with the vector bundle
$\pi^{\ast}:\varphi^{\ast}(\cA)\longrightarrow N$. The vector bundle $\pi_{TN}\circ \widehat{\gamma}:(T\varphi)^{\ast}\cA\longrightarrow N$
is $N$-isomorphic to a vector subbundle of the Whitney sum $TN\oplus\varphi^{\ast}(\cA)$ whose its sections are the sections
$$s=(X,\xi):N\longrightarrow TN\oplus\varphi^{\ast}(\cA)$$ ($X\in \Gamma(TN)$ and $\xi\in \Gamma(\varphi^{\ast}\cA)$) of the Whitney sum
$TN\oplus\varphi^{\ast}(\cA)$ characterized by the equality $T(\varphi)\circ X=\gamma\circ \Phi \circ \xi$, where $\Phi$ stands for the canonical
map from the vector bundle $\pi^{\ast}:\varphi^{\ast}(\cA)\longrightarrow N$ to the vector bundle $\pi:\cA\longrightarrow M$ defined by $\Phi(x,u)=u$.
The notable fact is that the vector bundle $\pi_{TN}\circ \widehat{\gamma}:(T\varphi)^{\ast}\cA\longrightarrow N$ inherits a natural structure of
transitive Lie algebroid on $N$. We note this structure of Lie algebroid on next proposition. The details of the proof can be found in
\cite{kuki-chw} or \cite{makz-lga}.

\vspace{3mm}

\begin{prop} (Inverse image of transitive Lie algebroids). If $\cA$ is a transitive
Lie algebroid on $M$, the vector bundle $\pi_{TN}\circ \widehat{\gamma}:(T\varphi)^{\ast}\cA\longrightarrow N$ carries a natural structure of
transitive Lie algebroid on $N$, in which the canonical projection $\widehat{\gamma}:(T\varphi)^{\ast}\cA\longrightarrow TN$ is the anchor
map and the Lie bracket is defined in the following way: Fix a local frame $(s_{1},\cdots,s_{k})$ of the vector bundle
$\pi:\cA\longrightarrow M$ defined on an open subset $U$ of $M$. Take now two sections $(X,\xi)$, $(Y,\eta)$ $\in \Gamma\big((T\varphi)^{\ast}\cA\big)$,
where $X,Y\in \Gamma(TN)$ and $\xi,\eta\in \Gamma(\varphi^{\ast}\cA)$. Then, on the open subset $V=\varphi^{-1}(U)$,
we have the decompositions $\xi_{/V}=\sum_{i}f_{i}(\s_{i}\circ \varphi_{/V})$ and $\eta_{/V}=\sum_{j}g_{j}(\s_{j}\circ \varphi_{/V})$ with
$f_{i},g_{j}\in C^{\infty}(V)$. Define a Lie bracket on $\Gamma\big((T\varphi)^{\ast}\cA\big)$ by setting
$$[(X,\xi),(Y,\eta)]_{/V}=$$
$$=\big([X,Y]_{/V},\sum_{i,j}f_{i}g_{j}[\s_{i},\s_{j}]\circ \varphi_{/V} +\sum_{j}(X\cdot g_{j})(\s_{j}\circ\varphi_{/V})
 -\sum_{i}(Y\cdot f_{i})(\s_{i}\circ\varphi_{/V})\big)$$ Thus, the pair of mappings $(\varphi^{!!},\varphi)$, in which
$\varphi^{!!}:(T\varphi)^{\ast}\cA\longrightarrow \cA$ is the smooth map defined by $\varphi^{!!}(X,a)=a$, is a morphism of Lie algebroids.
\end{prop}

\vspace{3mm}

Keeping the same hypothesis and notations as above, the vector bundle
$\pi_{TN}\circ \widehat{\gamma}:(T\varphi)^{\ast}\cA\longrightarrow N$, equipped with this structure of Lie algebroid, is called
the Lie algebroid inverse image of $\cA$ by the map $\varphi$ and denoted by $(\varphi^{!!}\cA,\widehat{\gamma},[\cdot,\cdot]_{\varphi^{!!}\cA})$.
From now on, we write simply $\varphi^{!!}\cA$ instead of $(\varphi^{!!}\cA,\widehat{\gamma},[\cdot,\cdot]_{\varphi^{!!}\cA})$,
dropping the anchor $\widehat{\gamma}$ and the Lie bracket $[\cdot,\cdot]_{\varphi^{!!}\cA}$. The map
$\varphi^{!!}:(T\varphi)^{\ast}\cA\longrightarrow \cA$ is called the canonical map induced and the Lie algebroid morphism $(\varphi^{!!},\varphi)$
is called the canonical Lie algebroid morphism of an induced Lie algebroid.

\vspace{3mm}

\textit{Example}. Let $\frak g$ be a finite dimensional Lie algebra. It is know that $\frak g$ is a Lie algebroid over a one-point space $M=\{\ast\}$. Let $N$ be a smooth manifold and $\varphi:N\longrightarrow M$ the constant map. Then, $(T\varphi)^{\ast}\frak g$ is equal to $TN\oplus (TN\times \frak g)$ and we easily see that the anchor and the Lie bracket of the Lie algebroid $(T(\varphi))^{\ast}\frak g$ coincides with the ones of the trivial Lie algebroid $TN\oplus (TN\times \frak g)$.

\vspace{3mm}

We define now the restriction of a transitive Lie algebroid to any submanifold of the base space.

\vspace{3mm}

\begin{defn} (Restriction of transitive Lie algebroids). Let $M$ be a smooth manifold and $\varphi:N\hookrightarrow M$ a submanifold, possibly with boundary and corners. Let $\cA$ be a transitive Lie algebroid on $M$. The Lie algebroid $\varphi^{!!}\cA$, constructed as inverse image of $\cA$ by the map $\varphi$, is called the Lie algebroid restriction of $\cA$ to the submanifold $N$ and denoted by $\cA^{!!}_{N}$.
\end{defn}

\vspace{3mm}

It is evident that, in the case in which $U$ is an open subset of $M$, the Lie algebroid restriction $\cA^{!!}_{U}$ has the natural structure of the transitive Lie algebroid given on the restricted vector bundle $\cA_{U}$, that is, the map $\psi:\cA_{U}\longrightarrow \cA^{!!}_{U}$ defined by $\psi(a)=(\gamma (a),a)$ is an $U$-isomorphism of Lie algebroids. Therefore, the Lie algebra structure on the set of the sections of $\cA^{!!}_{U}$ is defined by extending sections of $\cA^{!!}_{U}$ to sections of $\cA$. Next two propositions, we want to show that, if $M$ is a smooth manifold, $\varphi:N\hookrightarrow M$ a submanifold of $M$, which is a closed subset of $M$ in the topological sense, and $(\cA,[\cdot,\cdot],\gamma)$ a transitive Lie algebroid on $M$, the structure of Lie algebra on the set of the sections of $\cA^{!!}_{N}$ is also defined by natural extension of sections of $\cA^{!!}_{N}$ to sections of $\cA$.

\vspace{3mm}

\begin{prop} Let $M$ be a smooth manifold and $\varphi:N\hookrightarrow M$ a submanifold, possibly with boundary and corners.
Let $\cA$ be a transitive Lie algebroid on $M$, with anchor $\gamma:\cA\longrightarrow TM$ and Lie bracket $[\cdot,\cdot]$, and
$\pi:\cA\longrightarrow M$ the vector bundle underlying the Lie algebroid $\cA$. Denote by $\cA_{N}$ the vector bundle restriction of
$\cA$ to $N$ and $\varphi^{!!}:\cA^{!!}_{N}\longrightarrow \cA$ the canonical map defined by $\varphi^{!!}(X,a)=a$. Then, the following assertions holds.
\begin{enumerate}[a)]
\vspace{3mm}
\item The bundle $\im\varphi^{!!}$ is a vector subbundle of $\cA_{N}$ and
$\varphi^{!!}:\cA^{!!}_{N}\longrightarrow \im\varphi^{!!}$ is an $N$-isomorphism of vector bundles.
Its inverse is the map $(\varphi^{!!})^{-1}(a)=(\gamma (a),a)$.
\item $\im\varphi^{!!}=\gamma^{-1}(TN)$ and so the anchor $\gamma:\cA\longrightarrow TM$ restricts to $\im\varphi^{!!}\longrightarrow TN$.
\item If $\xi$, $\eta$ $\in \Gamma\cA$ are sections such that $\xi_{/N}$ and $\eta_{/N}$ $\in \Gamma(\im\varphi^{!!})$ then $[\xi,\eta]_{/N}$
is a section of $\im\varphi^{!!}$.
\item If $\xi$, $\eta$ $\in \Gamma\cA$ such that $\xi_{/N}$ $\in \Gamma(\im\varphi^{!!})$ and $\eta_{/N}=0$, then $[\xi,\eta]_{/N}=0$.
\end{enumerate}
\end{prop}

\vspace{2mm}

\proofwp. a) Standard arguments.

b) Let $a\in \im\varphi^{!!}$. Then, there exists $(x,u)\in TN$ such that $$\gamma(a)=(\varphi(x),D\varphi_{x}(u))=(x,u)\in TN$$ Hence, $\gamma(a)\in TN$.

c) Let $\xi,\eta\in \Gamma\cA$ such that $\xi_{/N}$, $\eta_{/N}$ $\in \Gamma(\im\varphi^{!!})$. Then, there are vector fields $X,Y\in TN$ such that
$(X,\xi_{/N})$, $(Y,\eta_{/N})$ $\in \Gamma(\cA^{!!}_{N})$, and so, $\gamma\circ \xi_{/N}=X$ and $\gamma\circ \eta_{/N}=Y$. Therefore, we have
$[(X,\xi_{/N}),(Y,\eta_{/N})]\in \Gamma(\cA^{!!}_{N})$. We will calculate a local expression for this last Lie bracket. Fix a local frame
$(s_{1},\cdots,s_{k})$ of the vector bundle $(\cA,\pi,M)$ over an open subset $U$ of $M$ and set $V=U\cap N$. Writing $\xi_{/U}=\sum_{i}f_{i}\s_{i}$
and $\eta_{/U}=\sum_{j}g_{j}\s_{j}$, with $f_{i},g_{j}\in C^{\infty}(U)$, we also can write $\xi_{/V}=\sum_{i}f_{i_{/V}}\s_{i/_{V}}$ and
$\eta_{/V}=\sum_{j}g_{j_{/V}}\s_{j/_{V}}$, and so, $[(X,\xi_{/N}),(Y,\eta_{/N})]_{/V}=\big([X,Y]_{/V},(\ast)\big)$ where $(\ast)$ means the following sum
$$(\ast)=\sum_{i,j}f_{i_{/V}}g_{j_{/V}}[\s_{i},\s_{j}]_{/V}+\sum_{j}(X_{/V}\cdot g_{j_{/V}})\s_{j_{/V}}\ -\sum_{i}(Y_{/V}\cdot f_{i_{/V}})\s_{i_{/V}}$$
On other side, $$[\xi,\eta]_{/U}=[\xi_{/U},\eta_{/U}]=[\xi_{/U},\sum_{j}g_{j}\s_{j}]=$$
$$=\sum_{j}g_{j}[\xi_{/U},s_{j}]+\sum_{j}\big((\gamma\circ \xi_{/U})\cdot g_{j})s_{j}\big)=$$
$$=-\sum_{j}g_{j}[s_{j},\xi_{/U}]+\sum_{j}\big(\sum_{i}f_{i}(\gamma\circ s_{i})\big)\cdot g_{j})s_{j}\big)=$$
$$=\sum_{i,j}f_{i}g_{j}[s_{i},s_{j}]-\sum_{j}g_{j}\big(\sum_{i}((\gamma\circ s_{j} )\cdot f_{i})s_{i}\big)+\sum_{j}\big(\sum_{i}f_{i}(\gamma\circ s_{i})\big)\cdot g_{j})s_{j}\big)=$$
$$=\sum_{i,j}f_{i}g_{j}[s_{i},s_{j}]-\sum_{i}\big(\sum_{j}(g_{j}(\gamma\circ s_{j}))\cdot f_{i})s_{i}\big)+\sum_{j}\big(\sum_{i}f_{i}(\gamma\circ s_{i})\big)\cdot g_{j})s_{j}\big)=$$ Since $X_{/V}=\gamma\circ \xi_{/V}=\gamma\circ (\sum_{i}f_{i_{/V}}\s_{i_{/V}})=\sum_{i}f_{i_{/V}}(\gamma\circ s_{i_{/N}})$ and, analogously \linebreak $Y_{/V}=\sum_{j}g_{j_{/V}}(\gamma\circ s_{j_{/N}})$ we conclude that $$[\xi,\eta]_{/V}=\sum_{i,j}f_{i_{/V}}g_{j_{/V}}[s_{i},s_{j}]_{/V}+$$ $$-\sum_{i}\big(\sum_{j}(g_{j_{/V}}(\gamma\circ s_{j_{/V}}))\cdot f_{i_{/V}})s_{i_{/V}}\big)+\sum_{j}\big(\sum_{i}f_{i_{/V}}(\gamma\circ s_{i_{/V}})\big)\cdot g_{j_{/V}})s_{j_{/V}}\big)=$$
$$=\sum_{i,j}f_{i_{/V}}g_{j_{/V}}[s_{i},s_{j}]_{/V}-\sum_{i}\big(Y_{/V})\cdot f_{i_{/V}})s_{i_{/V}}\big)+\sum_{j}\big(X_{/V}\cdot g_{j_{/V}})s_{j_{/V}}\big)=(\ast)$$ Therefore, $[(X,\xi_{/N}),(Y,\eta_{/N})]_{/V}=\big([X,Y]_{/V},[\xi,\eta]_{/V}\big)\in \Gamma((\cA^{!!}_{N})/_{V})$ and then $[\xi,\eta]_{/V}$ is a local section of the vector bundle $\im\varphi^{!!}$. Hence, $[\xi,\eta]$ is a section of the vector bundle $\im\varphi^{!!}$.

d) From part c) we have that $0=\gamma\circ \eta_{/N}=Y$ and so $[X,Y]=0$. From the equality $\eta_{/V}=\sum_{j}g_{j_{/V}}\s_{j/_{V}}$ we also have that
$g_{j_{/V}}=0$ for each $j$. Hence, on the expression $(\ast)$ of the part c), all summands are equal to zero. Then, $[\xi,\eta]_{/V}=0$ and so
$[\xi,\eta]_{/N}$=0 also. {\small $\square$}

\vspace{3mm}

\begin{prop} Keeping the same hypothesis and notations of the previous proposition, if $\cA$
is a transitive Lie algebroid on $M$ and $\varphi:N\hookrightarrow M$ a submanifold, which is a closed subset of $M$ in the topological sense,
then $\im\varphi^{!!}$ carries a natural structure of transitive Lie algebroid on $N$, where the anchor of this structure is the restriction
of $\gamma$ to $\im\varphi^{!!}\longrightarrow TN$ and the Lie bracket is defined by $$[\xi,\eta]=[\widehat{\xi},\widehat{\eta}]_{/N}$$ in
which $\widehat{\xi}$ and $\widehat{\eta}$ are sections of $\cA$ such that $\widehat{\xi}_{/N}=\xi$ and $\widehat{\eta}_{/N}=\eta$. Thus,
the map $\varphi^{!!}:\cA^{!!}_{N}\longrightarrow \im\varphi^{!!}$ is an $N$-isomorphism of Lie algebroids whose its inverse is the map
$\Psi:\im\varphi^{!!}\longrightarrow \cA^{!!}_{N}$ defined by $\Psi(a)=(\gamma(a),a)$.
\end{prop}

\proofwp. Since $N$ is a closed subset of $M$, it is always possible to extends sections of $\im\varphi^{!!}$ to a sections of $\cA$. By using c) and d) of last proposition and similar arguments given on the definition of Lie bracket of a restriction to an open subset, we can define a Lie bracket on the sections of $\im\varphi^{!!}$. {\small $\square$}

\vspace{3mm}

\textit{Remark}. The Lie algebroid $\im\varphi^{!!}$ constructed on last proposition can be identified to the Lie algebroid $\cA^{!!}_{N}$ and so, the Lie algebroid $\im\varphi^{!!}$ is also called the Lie algebroid restriction of $\cA$ to $N$. Henceforth, the Lie algebroid $\im\varphi^{!!}$ will also be denoted by $\cA^{!!}_{N}$.

\vspace{3mm}

We finalize this section with the following three propositions widely used in this paper.

\vspace{3mm}

\begin{prop} (Transitivity of restrictions) Let $M$, $N$ and $P$ three smooth manifolds such that $\varphi:N\hookrightarrow M$ is a submanifold,
which is a closed subset of $M$, and $\psi:P\hookrightarrow N$ is a submanifold, also closed in $N$ in the topological sense. Consider a transitive Lie
algebroid $(\cA,[\cdot,\cdot],\gamma)$ on $M$. Then, $$(\cA^{!!}_{N})^{!!}_{P}\simeq \cA^{!!}_{P}$$
\end{prop}

\proofwp. $$\cA^{!!}_{N}=\{(x,u,a)\in TN\times \cA: (x,u)=\gamma (a)\}$$ and $$(\cA^{!!}_{N})^{!!}_{P}=\{((y,v),(x,u),a)\in TP\times (TN\times \cA): (y,v)=\gamma_{\cA^{!!}_{N}}((x,u),a)\}$$ The map $\lambda$ from $\cA^{!!}_{P}$ onto $(\cA^{!!}_{N})^{!!}_{P}$ given by $\lambda(y,v,a)=((y,v),(y,v),a)$ is an isomorphism of Lie algebroids over $P$. {\small $\square$}

\vspace{3mm}

Given a smooth manifold $M$ and $N$ a submanifold of $M$, the vector bundle $TM_{N}$, restriction of $TM$ to a $N$, does not coincide in general with
the tangent bundle $TN$. However, in the context of restrictions of transitive Lie algebroids, the situation is different and much better because they coincide. We note that property in the next proposition.

\vspace{3mm}

\begin{prop} (Restriction of tangent Lie algebroids). Let $M$ be a smooth manifold and $\varphi:N\hookrightarrow M$ a compact embedded submanifold.
Consider the tangent Lie algebroids $TM$ and $TN$. Then, $(TM)^{!!}_{N}=TN$.
\end{prop}

\vspace{3mm}

Lastly, we have the following proposition.

\vspace{3mm}

\begin{prop} (Restriction of trivial Lie algebroids). Let $\frak g$ be a real finite dimensional Lie algebra and $M$ a smooth manifold.
Consider the trivial Lie algebroid $TM\oplus (M\times \frak g)$ and let $N$ be an embedded submanifold of $M$. Then,
$$(TM\oplus (M\times \frak g))_{N}^{!!}=TN\oplus (N\times \frak g)$$
\end{prop}

\proofwp. Let $\varphi:N\hookrightarrow M$ the inclusion. The Lie algebroid $(TM\oplus (M\times \frak g))_{N}^{!!}$, seen as the Lie algebroid $\im\varphi^{!!}$, is constituted by the elements $((x,u),(x,v))\in TM\oplus (M\times \frak g)$ such that there exists $(\widetilde{x},\widetilde{u})\in TN$ satisfying the equality $T(\varphi)(\widetilde{x},\widetilde{u})=\gamma((x,u),(x,v))$. Then, $(x,u)=(\widetilde{x},\widetilde{u})\in TN$ and so we have that $(TM \oplus (M\times \frak g))^{!!}_{N}=TN\oplus (N\times \frak g)$.
{\small $\square$}

\begin{center}
\section{Smooth forms and cohomology}
\end{center}

\vspace{6mm}

We shall recall briefly the notion of Lie algebroid cohomology and collect some results concerning Lie algebroid extensions. This section ends stating and proving the triviality of a transitive Lie algebroid on a contractible manifold. This last result will be used in the proof of the main theorem.

\vspace{3mm}

Let $M$ be a smooth manifold and $\cA$ a Lie algebroid on $M$. Let $\RR_{M}$ denote the trivial vector bundle $M\times \RR$
of base $M$. We recall that a smooth form of degree $p$ on $\cA$
is a section of $\bigwedge^{p}\big(\cA^{\ast};\RR_{M}\big)$. The set of
all smooth forms of degree $p$ on $\cA$ will be denoted by $\Omega^{p}(\cA;M)$. For $p=0$, we have
$\Omega^{0}(\cA;M)=C^{\infty}(M)$. The set $\Omega^{p}(\cA;M)$ is a $C^{\infty}(M)$-module for each $p\geq 0$. Clearly, the exterior product of alternated multi-linear maps
induces an exterior product in $\Omega^{\ast}(\cA;M)=\bigoplus_{p\geq 0} \Omega^{p}(\cA;M)$. Thus, this product makes $\Omega^{\ast}(\cA;M)$
into a commutative graded algebra in which the constant map
$\textrm{\textbf{1}}\in C^{\infty}(M)$ is the unit. We notice that $\Omega^{\ast}(\cA;M)$ vanishes
for degrees $>$ rank of $\cA$.

Next, we recall the definition of exterior derivative on $\cA$. We first consider the algebra $\Omega^{0}(\cA;M)=C^{\infty}(M)$. Let
$f\in C^{\infty}(M)$ a smooth map. We can define the smooth 1-form
$d(f):M\mapr{} \bigwedge^{1}\big(\Gamma
(TM),C^{\infty}(M)\big)$ by setting $d(f)(X)= X\cdot f$, for each $X\in
\Gamma (TM)$. Hence, we define $d(f)\in \bigwedge^{1}\big(\Gamma (\cA),C^{\infty}(M)\big)$ by
$d(f)(X)=(\gamma\circ X)\cdot f$, for each $X\in \Gamma (\cA)$. Now, for each $p\geq 1$ we define
$$d^{p}:\Omega^{p}(\cA;M)\longrightarrow \Omega^{p+1}(\cA;M)$$ $$d^{p}\omega
(X_{1},X_{2},\cdot\cdot\cdot,X_{p+1})=\sum_{j=1}^{p+1}(-1)^{j+1}(\gamma\circ
X_{j})\cdot(\omega
(X_{1},\cdot\cdot\cdot,\widehat{X_{j}},\cdot\cdot\cdot,X_{p+1}))\
+$$
$$+\sum_{i<k}(-1)^{i+k}\omega([X_{i},X_{k}],X_{1},\cdot\cdot\cdot,\widehat{X_{i}},\cdot\cdot\cdot,\widehat{X_{k}},\cdot\cdot\cdot,X_{p+1})$$
for $\omega \in \Omega^{p}(\cA;M)$ and
$X_{1},X_{2},\cdot\cdot\cdot,X_{p+1}\in \Gamma (\cA)$.

The family of differential operators $d^{\ast}=(d^{p})_{p\geq 0}$ defines, on the algebra $\Omega^{\ast}(\cA;M)$, a structure of differential
graded algebra, which is commutative. Hence, $\Omega^{\ast}(\cA;M)$ becomes a commutative cochain algebra defined over $\RR$.

\vspace{3mm}

\begin{defn} The Lie algebroid cohomology space of $\cA$ is the cohomology space of the cochain algebra $\Omega^{\ast}(\cA;M)$
equipped with the structures defined above. This cohomology space is denoted by $H^{\ast}(\cA;M)$.
\end{defn}

In order to state a result concerned with extensions of smooth forms on Lie algebroids, we recall the definition of inverse image of a smooth form. Let $\cA$ and $\cB$ be two Lie algebroids defined respectively over the smooth manifolds $M$ and $N$. Let $\lambda=(\psi,\varphi)$ be a morphism of Lie algebroids defined by the smooth maps $\psi:\cA\longrightarrow \cB$ and $\varphi:M\longrightarrow N$. If $\omega\in \Omega^{p} (\cB;N)$ is a smooth form on $\cB$ of degree $p$, we can consider a smooth form of degree $p$ on $\cA$, denoted by $\lambda^{\ast}\omega$, and defined by
$$\big(\lambda^{\ast}\omega\big)_{x}(v_{1},v_{2},\cdot\cdot\cdot,v_{p})=\omega_{\varphi(x)}(\psi(v_{1}),\psi(v_{2}),\cdot\cdot\cdot,\psi(v_{p}))$$
$x\in M$ and $v_{1}$, $v_{2}$, . . . $v_{p}$ $\in
\cA_{x}$. The form $\lambda^{\ast}\omega$ is called the pullback
or inverse image of $\omega$ by the morphism $\lambda$. Thus, for each $p\geq 0$, there is a map
$$\lambda^{\ast p}:\Omega^{p}(\cB;N)\longrightarrow \Omega^{p}(\cA;M)$$
$$\omega \longrightarrow \lambda^{\ast}\omega$$
The family $\lambda^{\ast}=(\lambda^{\ast p})_{p\geq 0}$ is a morphism of cochain algebras.

\vspace{3mm}

\begin{defn} (Restriction of smooth forms on Lie algebroids). Let $M$ be a smooth manifold and $\varphi:N\hookrightarrow M$ a submanifold.
Let $\cA$ be a Lie algebroid on $M$ and consider the canonical Lie algebroid morphism $\lambda=(\varphi^{!!},\varphi)$, in which
$\varphi^{!!}:\cA_{N}^{!!}\longrightarrow \cA$ is defined by $\varphi^{!!}(X,a)=a$. Consider the cochain algebras $\Omega^{\ast}(\cA;M)$
and $\Omega^{\ast}(\cA^{!!}_{N};N)$ and the morphism $\lambda^{\ast}:\Omega^{\ast}(\cA;M)\longrightarrow \Omega^{\ast}(\cA^{!!}_{N};N)$ induced.
For each smooth form $\omega\in \Omega^{\ast}(\cA;M)$, the form $\lambda^{\ast}(\omega)\in \Omega^{\ast}(\cA^{!!}_{N};N)$ is called the form
restriction of $\omega$ to $N$ and denoted by $\omega^{!!}_{N}$ or simply by $\omega_{/_{N}}$, if there is no danger of confusion with the
restriction of $\omega$ to the vector bundle restriction of the underlying vector bundle of $\cA$ to $N$. In subsequent sections,
the morphism $\lambda^{\ast}$ will often be denoted by $\varphi_{M,N}^{\cA,\cA^{!!}_{N}}$ or simply by $\varphi_{M,N}^{\cA}$.
\end{defn}

It is obvious that the restriction of smooth forms on Lie algebroids enjoys the functorial properties.

\vspace{3mm}

Let $\cA$ be a transitive Lie algebroid on a smooth manifold $M$. If $U$ an open subset of $M$ and $\varphi:U \longrightarrow M$ the inclusion map then, for each $p\geq 0$, the spaces $\Omega^{p}(\cA_{U};U)$ and $\Omega^{p}(\cA^{!!}_{U};U)$ are isomorphic. If $\varphi:N\hookrightarrow M$ a submanifold such that $N$ is a closed subset in $M$ in the topological sense, then the map $\varphi^{!!}:\cA^{!!}_{N}\longrightarrow \im\varphi^{!!}$ is a $N$-isomorphism of Lie algebroids and so it induces an isomorphism between the spaces $\Omega^{p}(\cA^{!!}_{N};N)$ and $\Omega^{p}(\im\varphi^{!!};N)$. Let us now notice a proposition concerning extensions of smooth forms.

\vspace{3mm}

\begin{prop} Let $M$ be a smooth manifold and $\varphi:N\hookrightarrow M$ a submanifold such that $N$ is a closed subset in $M$ in the
topological sense. Let $\cA$ be a transitive Lie algebroid on $M$ and consider the canonical Lie algebroid morphism $\lambda=(\varphi^{!!},\varphi)$
in which $\varphi^{!!}:\cA_{N}^{!!}\longrightarrow \cA$ is defined by $\varphi^{!!}(X,a)=a$. Then, the morphism of cochain algebras
$$\lambda^{\ast}:\Omega^{\ast}(\cA;M)\longrightarrow \Omega^{\ast}(\cA^{!!}_{N};N)$$ is surjective.
\end{prop}

\proofwp. Let us denote the anchor of $\cA$ by $\gamma$. Given a smooth form $\widetilde{\omega}\in \Omega^{p}(\cA^{!!}_{N};N)$,
we can define a smooth form $\widehat{\omega}\in \Omega^{p}(\im \varphi^{!!};N)$ by
$$\widehat{\omega}(\xi_{1},\dots,\xi_{p})=\widetilde{\omega}((\gamma\circ \xi_{1},\xi_{1}),\dots,(\gamma\circ \xi_{1},\xi_{p}))$$
The result follows since any smooth extension $\omega\in \Omega^{\ast}(\cA;M)$ of $\widehat{\omega}$ satisfies the equality
$\lambda^{\ast}(\omega)=\widetilde{\omega}$. {\small $\square$}

\vspace{3mm}

One of the properties listed in the introduction which is needed for the proof of our main result is the triviality of a transitive Lie algebroid over a contractible manifold. This property is a direct consequence of a deep result, due to Mackenzie, about actions from a certain cohomology space on the set of operators extensions of Lie algebroids by Lie algebra bundles. The notions of couplings corresponding extensions of Lie algebroids by Lie algebras bundles and the affine space of operator extensions play a leading role in the development of the matter in question. The heart of the theory is that the additive group of the cohomology space in degree two of a special representation induced by a coupling acts freely and transitively on the affine space of operator extensions
(see Mackenzie \cite{makz-lga}). In order to state both results, we briefly summarize now some definitions and examples used in the study of the subject concerned, following \cite{kuki-chw} and \cite{makz-lga}. Teleman presents in the paper \cite{Tele-rig} an approach to the theory of non-abelian Lie algebra extensions for the algebraic version of Lie algebroids. The MacLane's book \cite{mac-lane} contains an algebraic introduction to the theory of extensions.

\vspace{3mm}

Let $M$ be a smooth manifold and $\pi:E\longrightarrow M$ a vector bundle on $M$. For each $x\in M$, denote by $\cA(E)_{x}$ the vector space
of all linear maps $\psi:\Gamma (E)\longrightarrow E_{x}$ such that there exists a vector $u\in T_{x}M$ satisfying the equality
$$\psi(f\xi)=f(x)\psi(\xi)+(u\cdot f)_{x}\xi_{x}$$ for all $f\in C^{\infty}(M)$ and $\xi\in \Gamma(E)$. The vector $u$ is unique
and so we can define a map $$\gamma:\bigsqcup_{x\in M}\cA(E)_{x}\longrightarrow TM$$ We denote by $\cA(E)$ the disjoint union
$\bigsqcup_{x\in M}\cA(E)_{x}$. We define a Lie bracket on the space of the sections of $\cA(E)$ locally as follows.
Fix a local trivialization $\varphi:\pi^{-1}(U)\longrightarrow U\times F$ of the vector bundle $E$, in which $U$ is an open subset
of $M$ and $F$ is the fibre type of $E$. Let $\frak gl (F)$ be the Lie algebra of $F$ and, for each $\xi\in \Gamma (E)$, the map
$\xi_{\varphi}:U\longrightarrow F$ defined by $\xi_{\varphi}(x)=\varphi_{x}(\xi_{x})$, in which $\varphi_{x}:E_{x}\longrightarrow F$
is the linear map induced by $\varphi$. It can be seen in the section 1.2 of \cite{kuki-chw} that the map
$\overline{\varphi}:TU\times \frak gl (F)\longrightarrow \cA(E)_{U}$ defined by
$$\overline{\varphi}(u,g)(\xi)=\big(\varphi_{x}\big)^{-1}\big(u\cdot \xi_{\varphi} + (g\circ \xi_{\varphi})(x)\big)$$ is bijective.
Therefore, the anchor and the Lie bracket on the sections of the trivial Lie algebroid $TU\times \frak gl (F)$ can be carried to
$\cA(E)_{U}$ and the space $\cA(E)$ becomes a transitive Lie algebroid on $M$, which is denoted by $\cD(E)$ and called the
\textit{Lie algebroid of covariant differential operators on} $E$ (see \cite{kuki-chw}, \cite{makz-lga} and \cite{makz-gad}). Assume now that $\cA$ is a Lie algebroid on $M$, with anchor $\gamma:\cA\longrightarrow TM$.
A \textit{representation} $\rho$ of $\cA$ on $E$ is a morphism $\rho:\cA\longrightarrow \cD(E)$ of Lie algebroids. Under these conditions,
on the graded algebra $\Omega^{\ast}(\cA;M)$, an exterior derivative can be defined in the same way as we have defined exterior derivative above,
but taking $\rho (X_{j})$ instead of $\gamma\circ X_{j}$ (in our definition of exterior derivative given before, we used the trivial representation
$\rho'$ defined by $\rho'(X)=\gamma\circ X$). The cohomology space of this cochain algebra is denoted by $H(\cA,\rho,E)$.

\vspace{3mm}

\textit{Couplings}. Let $\cA$ be a Lie algebroid on a smooth manifold $M$ and $K$ a Lie algebra bundle on $M$. A coupling of $\cA$ with $K$
is a morphism $\Xi:\cA\longrightarrow \Out\cD(K)$ of Lie algebroids. In \cite{makz-lga} or \cite{makz-gad}, it can be seen that
a coupling induces a representation of $\cA$, denoted by $\rho^{\Xi}$, on the Lie algebra bundle $\Out\cD(ZK)$, in which $ZK$
denotes the Lie algebra subbundle center of $K$. The representation $\rho^{\Xi}$ is called central representation of the coupling $\Xi$.
A Lie derivation law covering the coupling $\Xi$ is a vector bundle morphism $\nabla:\cA\longrightarrow \cD_{Der}(K)$ that preserves
the anchor maps and satisfies the equality $\natural\circ \nabla=\Xi$, in which $\natural:\cD_{Der}(K)\longrightarrow \Out\cD(K)=\cD_{Der}(K)/\ad(K)$
is the quotient map. Since the map  $\natural$ is a surjective submersion, any coupling $\Xi$ admits Lie derivation laws covering it.
For transitive Lie algebroids, a Lie derivation law covering a couple $\Xi$ can be obtained by taking the covariant derivative of a
connection $\lambda:TM\longrightarrow \cA$.

\vspace{3mm}

\textit{Obstruction}. Let $\cA$ be a Lie algebroid on a smooth manifold $M$, $K$ a Lie algebra bundle on $M$ and
$\Xi:\cA\longrightarrow \Out\cD(K)$ a coupling of $\cA$ with $K$. For any Lie derivation law $\nabla:\cA\longrightarrow \cD_{Der}(K)$
covering the coupling $\Xi$, the curvature map $R_{\nabla}$ of $\nabla$ is the map
$R_{\nabla}:\Gamma(\cA)\times \Gamma(\cA)\longrightarrow \Gamma(\cD_{Der}(K))$ which is defined by
$$R_{\nabla}(\xi,\eta)=[\nabla(\xi),\nabla(\eta)]-\nabla([\xi,\eta])$$ for $\xi$, $\eta$ $\in \Gamma(\cA)$.
Let $\Lambda:\bigwedge^{2}(\cA)\longrightarrow K$ be a lifting of $R_{\nabla}$ (see \cite{mish-cla}, \cite{mish-latb} and \cite{mish-obs}). Then, the cyclic sum of $$\nabla_{\xi}(\Lambda(\eta,\theta))-\Lambda([\xi,\eta],\theta)$$ defines an element of $\Z^{3}(\cA,\rho^{\Xi},ZK)$.
The notable fact is that, the cohomology class of this element is independent of the choice of $\nabla$ and $\Lambda$, depending
only on the coupling $\Xi$. The cohomology class of this element is called the obstruction class of the coupling $\Xi$, and is
denoted by $Obs(\Xi)$ (see \cite{makz-lga}, theorem 7.2.12).

\vspace{3mm}

\textit{Extensions of Lie algebroids}. In view of the definition of set of equivalence classes of operator extensions, we recall that an exact sequence of Lie algebroids on a smooth manifold $M$ is a sequence
$$
\xymatrix{
\{0\} \ar[r] & \cA' \ar[r]^{j} & \cA \ar[r]^{\lambda} & \cA'' \ar[r] & \{0\}\\
}
$$ in which $\cA'$, $\cA$ and $\cA''$ are Lie algebroids on $M$, $j$ and $\lambda$ are morphisms of Lie algebroids and the sequence is exact as a sequence of vector bundles. We are interested in sequences in which $\cA'$ is a Lie algebra bundle. Given a Lie algebra bundle $K$, an extension of $\cA$ by $K$ is an exact sequence
$$
\xymatrix{
\{0\} \ar[r] & K \ar[r]^{j} & \cA' \ar[r]^{\lambda} & \cA \ar[r] & \{0\}\\
}
$$ of Lie algebroids over $M$. The most important extension for our work is the Atiyah sequence
$$
\xymatrix{
\{0\} \ar[r] & \ker \gamma \ar[r]^{j} & \cA \ar[r]^{\gamma} & TM \ar[r] & \{0\}\\
}
$$ of a transitive Lie algebroid $\cA$.

\vspace{3mm}

\textit{Equivalent extensions}. Let $M$ be a smooth manifold, $\cA$ a Lie algebroid on $M$ and $K$ a Lie algebra bundle on $M$. Two extensions
$$
\xymatrix{
\{0\} \ar[r] & K \ar[r]^{j_{1}} & \cA_{1} \ar[r]^{\lambda_{1}} & \cA \ar[r] & \{0\}\\
}
$$ and
$$
\xymatrix{
\{0\} \ar[r] & K \ar[r]^{j_{2}} & \cA_{2} \ar[r]^{\lambda_{2}} & \cA \ar[r] & \{0\}\\
}
$$ are equivalent if there is a Lie algebroid morphism $\varphi:\cA_{1}\longrightarrow \cA_{2}$ such that $\varphi\circ j_{1} =j_{2}$ and $\lambda_{2}\circ \varphi=\lambda_{1}$. In these conditions, $\varphi$ is an isomorphism of Lie algebroids.

\vspace{3mm}

\textit{Transversals}. Let $M$ be a smooth manifold, $\cA$ a Lie algebroid on $M$ and $K$ a Lie algebra bundle on $M$. Let
$$
\xymatrix{
\{0\} \ar[r] & K \ar[r]^{j} & \cA' \ar[r]^{\lambda} & \cA \ar[r] & \{0\}\\}$$ be an extension of $\cA$ by $K$. A transversal in the extension is a vector bundle morphism $\chi:\cA\longrightarrow \cA'$ such that $\lambda\circ \chi=id_{\cA}$. Since $\lambda$ is a surjective submersion and a morphism of vector bundles, transversals always exist and they are anchor-preserving morphisms. Fix now any transversal $\chi:\cA\longrightarrow \cA'$ in the extension above. We can define a coupling of $\cA$ with $K$ as follows. Define the map $\nabla^{\chi}:\cA\longrightarrow \cD_{Der}(K)$ such that $$j(\nabla^{\chi}(\xi)(Y))=[\chi(\xi),j(Y)]$$ for each $\xi\in \Gamma(\cA)$ and $Y\in \Gamma (K)$. The map $\nabla^{\chi}$ is a morphism of vector bundles preserving the anchor maps. The composition $\natural\circ \nabla^{\chi}: \cA\longrightarrow \Out\cD(K)$ is a coupling (see \cite{makz-lga}, section 7.3). If $\zeta:\cA\longrightarrow \cA'$ is other transversal of the extension, then we obtain the same coupling since the equality $\natural\circ \nabla^{\chi}=\natural\circ \nabla^{\zeta}$ holds. The coupling of $\cA$ with $K$ constructed in this way for any transversal is called the coupling induced by the extension. The map $\nabla^{\chi}$ is a Lie derivation law covering this coupling. Thus, in \cite{makz-lga}, it can be seen that, for any transversal in the extension and any Lie derivation law $\nabla$ covering the coupling of $\cA$ with $K$ induced by the extension, there is other transversal $\zeta:\cA\longrightarrow \cA'$ such that $\nabla=\nabla^{\zeta}$.

\vspace{3mm}

\textit{Operator extensions}. Let $M$ be a smooth manifold, $\cA$ a Lie algebroid on $M$ and $K$ a Lie algebra bundle on $M$. Let  $\Xi: \cA\longrightarrow \Out\cD(K)$ be a coupling of $\cA$ with $K$ such that $Obs(\Xi)=0\in \H^{3}(\cA,\rho^{\Xi},ZK)$, in which $\rho^{\Xi}$ denotes the central representation of the coupling. An operator extension of $\cA$ with $K$ is an extension
$$
\xymatrix{
\{0\} \ar[r] & K \ar[r]^{j} & \cA' \ar[r]^{\lambda} & \cA \ar[r] & \{0\}\\
}
$$ such that the coupling induced by this extension coincide with the coupling  $\Xi$. The set of equivalence classes of operator extensions of $\cA$ by $K$ is denoted by $\cO(\cA,\Xi,K)$.

\vspace{3mm}

The result, due to Mackenzie, is that $\H^{2}(\cA,\rho^{\Xi},ZK)$ acts freely and transitively on the set $\cO(\cA,\Xi,K)$. We begin by defining this action.

\vspace{3mm}

\textit{Action on operation extensions}. Let $M$ be a smooth manifold, $\cA$ a Lie algebroid on $M$ and $K$ a Lie algebra bundle on $M$. Let $\Xi: \cA\longrightarrow \Out\cD(K)$ be a coupling of $\cA$ with $K$ such that $Obs(\Xi)=0\in \H^{3}(\cA,\rho^{\Xi},ZK)$, in which $\rho^{\Xi}$ denotes the central representation of the coupling. Consider an operator extension
$$
\xymatrix{
\{0\} \ar[r] & K \ar[r]^{j} & \cA' \ar[r]^{\lambda} & \cA \ar[r] & \{0\}\\
}
$$ Let $g\in \Z^{2}(\cA;ZK)$. Then, the action of $g$ on the extension yields the extension
$$
\xymatrix{
\{0\} \ar[r] & K \ar[r]^{j} & \cA'_{g} \ar[r]^{\lambda} & \cA \ar[r] & \{0\}\\
}
$$ in which $\cA'_{g}=\cA'$ as vector bundles, the maps $j$ and $\lambda$ are the same in both extensions, the anchors $\gamma':\cA'\longrightarrow TM$ and $\gamma'_{g}:\cA'_{g}\longrightarrow TM$ are the same too and the Lie bracket $[\cdot,\cdot]_{g}$ on $\Gamma(\cA'_{g})$ is given by $$[\xi,\eta]_{g}=[\xi,\eta]+(j\circ i\circ g)(\lambda(\xi),\lambda(\eta)$$ in which $i:ZK\longrightarrow K$ denotes the inclusion.

\vspace{3mm}

We can now state the main result concerning actions of operators extensions. All details of the (long) proof can be found in Mackenzie \cite{makz-lga}, section 7.3.

\vspace{3mm}

\begin{prop} (Mackenzie). Let $\cA$ be a Lie algebroid on a smooth manifold $M$, $K$ a Lie algebra bundle on $M$ and
$\Xi: \cA\longrightarrow \Out\cD(K)$ a coupling of $\cA$ with $K$. Denote by $\rho^{\Xi}$ the central representation corresponding
the coupling $\Xi$. Suppose that the class obstruction of $\Xi$ is the zero of $\H^{3}(\cA,\rho^{\Xi}, ZK)$. Then,
the additive group of $\H^{2}(\cA,\rho^{\Xi}, ZK)$ acts freely and transitively on $\cO(\cA,\Xi,K)$.
\end{prop}

\vspace{3mm}

There are several consequences from the previous proposition. For our work, the most important consequence is the triviality of a transitive
Lie algebroid on a contractible manifold.

\vspace{3mm}

\begin{prop} Let $\cA$ be a transitive Lie algebroid on a contractible smooth manifold $M$ and $\frak g=\ker \gamma$,
in which $\gamma$ denotes the anchor of $\cA$. Then, $\cA$ is isomorphic to the trivial Lie algebroid $TM\times \frak g$ by a
strong isomorphism of Lie algebroids.
\end{prop}

\proofwp. Since $\cA$ is transitive, we can fix a connection $a:TM\longrightarrow \cA$ (see \cite{makz-lga} and \cite{makz-gad}). Let $\Xi: \cA\longrightarrow \Out\cD(K)$ be the coupling defined by $\Xi=\natural \circ \nabla^{a}:\cA\longrightarrow \Out\cD(K)$, in which $\natural$ is the quotient map $\natural: \frak D _{Der}\frak g \longrightarrow \Out\cD(K)$. Consider the extensions
$$
\xymatrix{
\{0\} \ar[r] & K \ar[r]^{j} & \cA \ar[r]^{\gamma} & TM \ar[r] & \{0\}
}
$$
and
$$
\xymatrix{
\{0\} \ar[r] & K \ar[r]^(.3){j} & TM\times \frak g \ar[r]^(.6){\gamma} & TM \ar[r] & \{0\}
}
$$
Since $M$ is contractible then $\H^{3}(TM,\rho^{\Xi}, ZK)=\{0\}$ and, by transitivity of the action of the previous proposition, there exists $[g]\in \H^{2}(TM,\rho^{\Xi}, ZK)$ such that the extensions
$$
\xymatrix{
\{0\} \ar[r] & \frak g \ar[r]^{j} & \cA_{g} \ar[r]^{\gamma} & TM \ar[r] & \{0\}
}
$$
and
$$
\xymatrix{
\{0\} \ar[r] & \frak g \ar[r]^(.3){j} & TM\times \frak g \ar[r]^(.6){\gamma} & TM \ar[r] & \{0\}
}
$$ are equivalent. Since $M$ is contractible, we can take $h\in \Omega^{1}(TM,\rho^{\Xi}, ZK)$ such that $dh=g$. It is known that the extensions
$$
\xymatrix{
\{0\} \ar[r] & \frak g \ar[r]^{j} & \cA_{dh} \ar[r]^{\gamma} & TM \ar[r] & \{0\}
}
$$
and
$$
\xymatrix{
\{0\} \ar[r] & \frak g \ar[r]^{j} & \cA_{g} \ar[r]^{\gamma} & TM \ar[r] & \{0\}
}
$$ are equivalent (see (\cite{makz-lga}). The conclusion follows by transitivity. {\small $\square$}

\vspace{3mm}

A proof of this result has been done by Mackenzie in \cite{makz-lga}, following the theory of non-abelian Lie algebra extensions. Crainic and Fernandes proved the same result (see \cite{crac-loja}, corollary 5.6) by using their theorem on integrability of Lie algebroids (see \cite{crac-loja}, corollary 5.4 and theorem 4.1).

\begin{center}

\item\section{Piecewise smooth forms and cohomology}

\end{center}

\vspace{6mm}

For each simplicial complex $K$, its polytope will be denoted by $|K|$. Simplex means always closed simplex. A smooth manifold $M$ is said to be smoothly triangulated by a simplicial complex $K$ if there is a homeomorphism $t$ from $|K|$ onto $M$ and, for each simplex $\Delta$ of $K$, the map $t_{/\Delta}:\Delta\longrightarrow M$ is a smooth embedding. As usually is done, we shall not make no notational distinction between the manifold and the complex triangulating it. In what follows, all simplicial complexes considered are geometric and finite. Each simplex can be represented as the convex body generated
by its vertices and, if its vertices are the points $a_{0}$, $a_{1}$, $\dots$, $a_{p}$, we eventually denote this simplex by
$[a_{0},a_{1},\dots,a_{p}]$. We shall write $\Delta' \prec\Delta$, if $\Delta'$ is a face of the simplex $\Delta$. The notation
$\varphi:\Delta'\hookrightarrow \Delta$, where $\varphi$ is the inclusion map, also will be used when $\Delta'$ is a face of $\Delta$.
The star of the simplex $\Delta$ in a simplicial complex $K$, denoted $\St(\Delta)$, is the union of the interiors of all (closed)
simplices of $K$ having $\Delta$ as a face. The closed star of the simplex $\Delta$ in $K$ is the union of all (closed) simplices of
$K$ having $\Delta$ as a face. The star $\St (\Delta)$ is an open subset in $|K|$ for the weak topology (which is the same as the
topology of subspace induced by the topology of the ambient space since the simplicial complex is finite).

There are several variants of piecewise smooth cohomology. What we consider here is mainly the piecewise smooth cohomology of transitive Lie algebroids
defined over smooth manifolds which are smoothly triangulated by simplicial complexes. It is in this context that our main result arises.
Its proof goes beyond over this notion of piecewise smooth cohomology. Precisely, in the course of the proof we will deal with the piecewise
smooth cohomology of Lie algebroids defined over bases which are unions of open stars in the polytope of the simplicial complex. Unions of stars were called by Eilenberg-Steenrod regular neighborhoods and we will follow this terminology (for general definition of regular neighborhood, see Eilenberg-Steenrod \cite{Eile-stin}, section 9 of the second chapter).

\vspace{3mm}

Let $M$ be a smooth manifold, smoothly triangulated by a simplicial complex $K$, and $\cA$ a transitive Lie algebroid on $M$. Let $\Delta$ be a simplex of $K$. Since $\Delta$ is an embedded submanifold of $M$, the Lie algebroid restriction $\cA_{\Delta}^{!!}$ is well defined. Suppose that $\Delta'$ is other simplex of $K$ such that $\varphi_{\Delta,\Delta'}:\Delta'\hookrightarrow \Delta$ is a face of $\Delta$. By transitivity of restrictions, we have that $\cA_{\Delta'}^{!!}\simeq(\cA_{\Delta}^{!!})^{!!}_{\Delta'}$. Consequently, the cochain algebras $\Omega^{\ast}(\cA_{\Delta'}^{!!};\Delta')$ and $\Omega^{\ast}((\cA_{\Delta}^{!!})^{!!}_{\Delta'};\Delta')$ are isomorphic.  We recall that the morphism of cochain algebras generated by the inclusion $\varphi_{\Delta,\Delta'}:\Delta'\hookrightarrow \Delta$ is denoted by $$\varphi^{\cA_{\Delta}^{!!}}_{\Delta,\Delta'}:\Omega^{\ast}(\cA_{\Delta}^{!!};\Delta)\longrightarrow \Omega^{\ast}(\cA_{\Delta'}^{!!};\Delta')$$ and, for each smooth form $\omega_{\Delta}\in \Omega^{p}(\cA_{\Delta}^{!!};\Delta)$, the smooth form $\varphi^{\cA_{\Delta}^{!!}}_{\Delta,\Delta'}(\omega_{\Delta})$ is also denoted by $(\omega_{\Delta})^{!!}_{\Delta'}$ or $(\omega_{\Delta})_{/_{\Delta'}}$. Keeping this hypothesis and notations, we give below the definition of piecewise smooth form. The idea of this definition is based in the Whitney's book \cite{wity-git} or in the Sullivan's papers \cite{suli-inf} and \cite{suli-Tok}. Morgan and Griffiths has also presented in \cite{gang-diff} the notion of piecewise smooth form on an ambient space made by a collection of manifolds with transverse intersections.

\vspace{3mm}

\begin{defn} (Piecewise smooth form). Let $M$ be a smooth manifold, smoothly triangulated by a simplicial complex $K$, and $\cA$ a transitive Lie algebroid on $M$. A piecewise smooth form of degree $p$ ($p\geq0$) on $\cA$ is a family $\omega=(\omega_{\Delta})_{\Delta \in K}$ such that the
following conditions are satisfied.

\begin{itemize}
\item For each $\Delta\in K$, $\omega_{\Delta}\in \Omega^{p}(\cA_{\Delta}^{!!};\Delta)$ is a smooth
form of degree $p$ on $\cA_{\Delta}^{!!}$.
\item For each $\Delta$, $\Delta'$ $\in K$, if $\varphi_{\Delta,\Delta'}:\Delta'\hookrightarrow \Delta$ is a face of $\Delta$, $$\varphi^{\cA_{\Delta}^{!!}}_{\Delta,\Delta'}(\omega_{\Delta})=\omega_{\Delta'}$$
\end{itemize}
\end{defn}

\vspace{3mm}

Let $(\varphi_{\Delta,\Delta'})^{!!}:\cA_{\Delta'}^{!!}\longrightarrow \cA_{\Delta}^{!!}$ be the map induced by $\varphi_{\Delta,\Delta'}$ (see the paragraph after proposition 1.1). The Lie algebroid $\cA_{\Delta'}^{!!}$ can be identified to the Lie algebroid $\im (\varphi_{\Delta,\Delta'})^{!!}$ (see remark after proposition 1.4) and so, for each $x\in \Delta'$, the fibre $(\cA_{\Delta'}^{!!})_{x}$ is a vector subspace of the fibre $(\cA_{\Delta}^{!!})_{x}$. Besides, the spaces $\Omega^{\ast}(\cA_{\Delta'}^{!!};\Delta')$ and $\Omega^{\ast}(\im (\varphi_{\Delta,\Delta'})^{!!};\Delta')$ are identified. Hence, the second condition of the definition given above can be stated in the following form: for each
$x\in \Delta'$ and vectors $u_{1}$, . . . , $u_{p}$ $\in
(A_{\Delta'}^{!!})_{x}$ $$\omega_{\Delta'}(x)(u_{1}, \cdot\cdot\cdot,
u_{p})=\omega_{\Delta}(x)(u_{1}, \cdot\cdot\cdot, u_{p})$$
Thus, a piecewise smooth form on $\cA$ is a collection of smooth forms, each one defined on the Lie algebroid restriction of $\cA$ to each simplex of $K$, which are compatible under restriction to faces. The set of all piecewise smooth forms of
degree $p$ on $\cA$ will be denoted by $\Omega^{p}(\cA;K)$. When we need to emphasize the piecewise context, we write $\Omega^{p}_{\textrm{ps}}(\cA;M)$ instead of $\Omega^{p}(\cA;K)$. We have then
$$\Omega^{p}(\cA;K)=\{(\omega_{\Delta})_{\Delta\in
K}:\omega_{\Delta}\in \Omega^{p}(\cA^{!!}_{\Delta};\Delta), \ \ \Delta'\prec\Delta\
\ \Longrightarrow \ \ (\omega_{\Delta})_{/_{\Delta'}}=\omega_{\Delta'}\}$$

When $p=0$, a piecewise smooth form of degree zero on $\cA$
is a family $(\varphi_{\Delta})_{\Delta\in K}\in \prod_{\Delta\in K}C^{\infty}(\Delta)$ such that
$\varphi_{\Delta}:\Delta\longrightarrow  \RR$ is smooth and the equality $\varphi_{\Delta'}=\varphi_{\Delta_{/\Delta'}}$ holds for each face $\Delta'$ of $\Delta$. The compatibility condition of restrictions to faces gives a map $\varphi:|K|\longrightarrow \RR$ which is continuous. The map $\varphi$ may not be a differentiable map but it is a piecewise smooth function. The set of all maps $\varphi\in C(|K|;\RR)$ which are compatible with restrictions to the faces of $|K|$ and with smooth restrictions to the faces of $|K|$ is denoted by $C_{ps}(|K|;\RR)$. Obviously, $\Omega^{0}(\cA;K)$ has a natural structure of algebra over $\RR$ and is naturally identified to $C_{ps}(|K|;\RR)$.

Since the restrictions of smooth forms are compatible with sums and products, various
operations on $\Omega^{p}(\cA;K)$ can be defined by the corresponding operations on each simplex of $K$. The set $\Omega^{p}(\cA;K)$, equipped with these operations, becomes a real vector subspace of $\prod_{\Delta\in K}\Omega^{p}(\cA_{\Delta};\Delta)$, for each natural $p\geq 0$. Thus, $\Omega^{p}(\cA;K)$ is a module over the algebra $C_{ps}(|K|;\RR)$. When $p=0$, $\Omega^{0}(\cA;K)=C_{ps}(|K|;\RR)$ has a structure of an unitary associative algebra over $\RR$. Moreover, the direct sum $$\Omega^{\ast}(\cA;K)=\bigoplus_{p\geq 0}\Omega^{p}(\cA;K)$$ equipped with the exterior product defined by the corresponding exterior product on each algebra $\Omega^{\ast}(\cA_{\Delta};\Delta)=\bigoplus_{p\geq 0}\Omega^{p}(\cA_{\Delta};\Delta)$, is a commutative graded algebra over $\RR$.

In order to obtain a complex of cochains, especially important is the analogues of exterior derivative. This operator also is obtained by the corresponding exterior derivative on the restriction of $\cA$ to each simplex of $K$. Such as in the case of smooth forms on a Lie algebroid, the space $\Omega^{\ast}(\cA;K)$, with the operations and differentiation above, becomes a cochain algebra, which is defined over $\RR$.

\vspace{3mm}

\begin{defn} (Piecewise smooth cohomology). Keeping the same hypothesis and notation as above, the piecewise smooth cohomology space of $\cA$ is the cohomology space of the cochain algebra $\Omega^{\ast}(\cA;K)$ equipped with the structures defined above. Its cohomology, $H(\Omega^{\ast}(\cA;K))$, will be denoted by $H^{\ast}(\cA;K)$ or $H^{\ast}_{ps}(\cA;M)$.
\end{defn}

\vspace{3mm}

We shall formulate now the main problem of this paper. Let $M$ be a compact smooth manifold, smoothly triangulated by a simplicial complex $K$ and $\cA$ a transitive Lie algebroid on $M$. Let $\omega \in \Omega^{p}(\cA; M)$ be a smooth form on $\cA$ of degree $p$. For each simplex $\Delta$ of $K$, let $\varphi_{M,\Delta}:\Delta\longrightarrow M$ be the inclusion map. We can restrict the form $\omega$ to the smooth form $\omega_{/_{\Delta}}=\varphi^{\cA}_{M,\Delta}(\omega)\in \Omega^{\ast}(\cA^{!!}_{\Delta};\Delta)$. It is obvious that the family $\omega=(\omega_{/_{\Delta}})_{\Delta \in K}$ is a piecewise smooth form on $\cA$. Hence, we have a linear map $$\Psi^{p}:\Omega^{p}(\cA;M)\longrightarrow \Omega^{p}(\cA;K)$$
defined by $$\omega \longrightarrow (\omega_{/_{\Delta}})_{\Delta\in
K}$$ Since the exterior derivative commutes with restrictions to any submanifold of $M$, the family $\Psi=(\Psi^{p})_{p\geq 0}$ defines a morphism of cochain algebras from $\Omega^{\ast}(\cA;M)$ to $\Omega^{\ast}(\cA;K)$. This map $\Psi$ is called the restriction map. We claim that \textbf{the restriction map $\Psi$ induces an isomorphism in cohomology}. This is our main result and the paper is devoted to the proof of this result.

\vspace{3mm}

The relationship between Lie algebroid and piecewise smooth cohomology of a transitive Lie algebroid on a triangulated compact manifold will be made through this isomorphism. There are some facts that we will need for the proof of this theorem such that the Mayer-Vietoris sequences for regular open subsets in the smooth and piecewise smooth contexts, the K$\ddot{\textrm{u}}$nneth theorem in both contexts, the triviality of Lie algebroids over contractible manifolds and the de Rham-Sullivan theorem for cell manifolds. It should be remarked that, for the statement of the Mayer-Vietoris sequence in the piecewise smooth context, we are going to deal with a complex of piecewise smooth forms which may not be defined over a family of closed simplices of a simplicial complex but over a collection of submanifolds obtained by the intersection of closed simplices with a regular open subset. Note that, the definitions and constructions we have made within the context of simplicial complexes did not require specific properties of the simplices. Indeed, these notions and properties can be extended to more general context. To illustrate this idea, let us briefly look at some cases of these constructions. A first example of this generalization is take a simplicial complex and to fix our attention on an open star of one its vertex. The family of submanifolds made by those simplices without the faces opposite to the vertex will be used for our generalization of the piecewise smooth setting. Another illustrative example consists of taking the family defined by intersections of open stars with any open subset of the polytope of a simplicial complex. In this case, the construction of the piecewise smooth setting is done in similar way. The first example is obviously a particular case of this second example. These two examples will be used in the proof of the main result.

We restrict now our attention to this generalization of the piecewise smooth setting. We begin by providing some definitions and notations. Once these ideas are established, we shall then turn towards to the statement of the Mayer-Vietoris sequences for regular open subsets in the smooth and piecewise smooth cases.

\vspace{3mm}

\textit{Generalization of the piecewise setting}. Let $M$ be a smooth manifold, smoothly triangulated by a simplicial complex $K$. Assume that $\s_{1}$, $\dots$, $\s_{e}$ are simplices of $K$ and consider the open subsets
$$U_{1}=\St(\s_{1}),\dots, U_{e}=\St(\s_{e})$$
Let $U=U_{1}\cup \dots \cup U_{e}$. We can consider the family $\underline{K}^{U}$ of all submanifolds $\Delta_{U}=\Delta\cap U$ such that
$\Delta$ is a simplex of $K$ provide that $\s_{j}$ is a face of $\Delta$ for some $j\in \{1,\dots,e\}$, that is,
$$\underline{K}^{U}=\big\{U\cap \Delta: \Delta\in K, \s_{j}\prec \Delta, j\in \{1,\dots,e\}\big\}$$
Suppose that $\cA$ is a transitive Lie algebroid on $U$. A piecewise smooth form of degree $p$ on $\cA$ is a family $$\omega=(\omega_{\Delta_{U}})_{\Delta_{U}\in \underline{K}^{U}}\in\prod _{\Delta_{U}\in \underline{K}^{U}}\Omega^{p}(\cA_{\Delta_{U}}^{!!};\Delta_{U})$$ such that, if $\Delta$ and $\Delta'$ are simplices of $K$,
with $\s_{j}\prec \Delta'\prec \Delta$ for some $j\in \{1,\dots,e\}$, and $\varphi_{\Delta_{U},\Delta'_{U}}:\Delta'_{U}\longrightarrow \Delta_{U}$ the inclusion map, one has
$$\varphi^{\cA^{!!}_{\Delta_{U}}}_{\Delta_{U},\Delta'_{U}}(\omega_{\Delta_{U}})=\omega_{\Delta_{U}'}$$ or simply $(\omega_{\Delta_{U}})_{/\Delta_{U}'}=\omega_{\Delta_{U}'}$.

\vspace{3mm}

The set $\Omega^{\ast}_{ps}(\cA;U)$ of all piecewise smooth forms on $\cA$ is a graded real vector space. A wedge product and a differential can be defined on $\Omega^{\ast}_{ps}(\cA;U)$ by the corresponding operations on each cochain algebra $\Omega^{\ast}(\cA_{\Delta_{U}}^{!!};\Delta_{U})$, giving to $\Omega^{\ast}_{ps}(\cA;U)$ a structure of cochain algebra defined over $\RR$. The cohomology space of this cochain algebra is denoted by $H^{\ast}_{ps}(\cA;U)$.
As before, we can define a restriction map $$\Psi:\Omega^{\ast}(\cA;U)\longrightarrow \Omega^{\ast}_{ps}(\cA;U)$$
by $$\omega \longrightarrow (\omega_{/_{\Delta_{U}}})_{\Delta_{U}\in
\underline{K}^{U}}$$ This map $\Psi$ is a morphism of cochain algebras.

\vspace{3mm}

Next, we are concerned with a piecewise variant of the Mayer-Vietoris sequence. As in the smooth case, the Mayer-Vietoris sequence is the long sequence induced from the canonical short exact sequence corresponding to two open subsets with union equal to the space and maps given by restriction and difference of forms. We begin by stating this short exact sequence in the piecewise context.

\vspace{3mm}

Let $M$ be a smooth manifold, smoothly triangulated by a simplicial complex $K$, and $\cA$ a transitive Lie algebroid on $M$. Let $\s_{1}$, $\dots$, $\s_{e}$ be simplices of $K$ and consider the open subsets $U_{j}=\St(\s_{j})$, with $j\in \{1,\dots,e\}$. For $l\in \{1, \dots, e\}$
fixed, consider the open subsets $U=U_{1}\cup \dots \cup U_{l}$ and
$V=U_{l+1}\cup \dots \cup U_{e}$ of $M$ and assume that $M=U\cup V$. Consider the following sets of manifolds:
\begin{itemize}
\item $\underline{K}^{U}$ is the set of all submanifolds $\Delta_{U}=U\cap \Delta$ such that
$\Delta\in K$ and $\s_{j}$ is a face of $\Delta$ for some $j\in \{1, \dots, l\}$;
\item $\underline{K}^{V}$ is the set of all submanifolds $\Delta_{V}=V\cap \Delta$ such that
$\Delta\in K$ and $\s_{i}$ is a face of $\Delta$ for some $i\in \{l+1, \dots, e\}$;
\item $\underline{K}^{U\cap V}$ is the set of all submanifolds $\Delta_{U\cap V}=(U\cap V)\cap \Delta$ such that
$\Delta\in K$ and $\s_{j}$ and $\s_{i}$ are faces of $\Delta$ for some $j\in \{1, \dots, l\}$ and $i\in \{l+1, \dots, e\}$.
\end{itemize}
Define the following two maps $\delta$ and $\pi$,
$$\delta:\Omega^{p}(\cA;K)\longrightarrow
\Omega^{p}_{ps}(\cA_{U};U)\times
 \Omega^{p}_{ps}(\cA_{V};V)$$
$$\pi:\Omega^{p}_{ps}(\cA_{U};U)\times
\Omega^{p}_{ps}(\cA_{V};V)\longrightarrow
\Omega^{p}_{ps}(\cA_{U\cap V};U\cap V)$$
by
$$\delta \Big((\omega_{\Delta})_{\Delta\in K}\Big)=
\Big(\big(\omega_{\Delta_{/\Delta_{U}}}\big)_{\Delta_{U}\in \underline{K}^{U}},\big(\omega_{\Delta_{/\Delta_{V}}}\big)_{\Delta_{V}\in \underline{K}^{V}}\Big)$$
$$\pi
\Big((\xi_{\Delta_{U}})_{\Delta_{U}\in \underline{K}^{U}},(\eta_{\Delta_{V}})_{\Delta_{V}\in \underline{K}^{V}}\Big)=\Big((\eta_{\Delta_{V}})_{/\Delta_{U\cap V}}-(\xi_{\Delta_{U}})_{/\Delta_{U\cap V}})_{\Delta_{U\cap V}\in \underline{K}^{U\cap V}}\Big)$$

\vspace{3mm}

Keeping these hypothesis and notation, we present our next result.

\vspace{3mm}

\begin{prop} The sequence
\begin{equation*}
\begin{array}{ccccc}
\{0\}\longrightarrow \Omega^{p}(\cA;K) & \mapr{\delta} &
\Omega^{p}_{ps}(\cA_{U};U)\oplus \Omega^{p}_{ps}(\cA_{V};V) & \mapr{\pi} & \Omega^{p}_{ps}(\cA_{U\cap V};U\cap V)\longrightarrow \{0\}\\
\end{array}
\end{equation*}
is a short exact sequence.
\end{prop}

\proofwp. Let $\Delta$ be a simplex of $K$ and $v$ the barycenter of this simplex. Since $M=U\cup V$, the collection of the stars $\{U_{j}: j\in \{1,\dots,e\}\}$ is an open covering of $M$. Then, there is an index $j\in \{1,\dots,e\}$ and a simplex $\Delta'$ of $K$ such that $s_{j}$ is a face of $\Delta'$ and $v$ belongs to the interior of $\Delta'$. Hence, the simplex $s_{j}$ is a face of $\Delta$. Therefore, each simplex of $K$ has some $s_{j}$ as a face. With this property, we easily can see that the map $\delta$ is injective and that $\im \delta=\ker \pi$. We shall check now that the map $\pi$ is surjective. Since
the set $\{U,V\}$ is an open covering of $M$, we can fix two smooth
maps $\varphi,\psi:M\longrightarrow [0,1]$ such that
$\supp \varphi\subset U$, $\supp \psi\subset V$ and
$\varphi(x)+\psi(x)=1$ for each $x$ in $M$. Let
$$(\gamma_{\Delta_{U\cap V}})_{\Delta_{U\cap V}\in \underline{K}^{U\cap V}}\in
\Omega^{p}_{ps}(\cA_{U\cap V};U\cap V)$$ be a piecewise smooth form on $\cA_{U\cap V}$. We shall define a
differential form $$(\xi_{\Delta_{U}})_{\Delta_{U}\in \underline{K}^{U}}\in \Omega^{p}_{ps}(\cA_{U};U)$$ as follows. For each $\Delta_{U}\in \underline{K}^{U}$, set
$$\xi_{\Delta_{U}}(x)=\begin{cases}

-\psi(x)\ \gamma_{\Delta_{U\cap V}}(x) & \text{if $x\in
\Delta\cap U\cap V$}
\\\\\
\ 0_{x}\in (\cA^{!!}_{\Delta_{U}})_{x} & \text{if $x\in
\Delta\cap U\cap(M\setminus \supp \psi)$}
\end{cases}$$
The sets $\Delta_{U}\cap V$ and
$\Delta_{U}\cap(M\setminus \supp \psi)$ are open in
$\Delta_{U}$ with union equal to $\Delta_{U}$.
Obviously, the restrictions of $\xi_{\Delta_{U}}$ to
$\Delta_{U}\cap V$ and to $\Delta_{U}\cap(M\setminus \supp \psi)$ are smooth. Therefore, we conclude that
$\xi_{\Delta_{U}}\in \Omega^{p}(\cA^{!!}_{\Delta_{U}};\Delta_{U})$.
In order to obtain a piecewise smooth form belonging to $\Omega^{p}_{ps}(\cA_{U};U)$, it remains to check that
$(\xi_{\Delta_{U}})_{\Delta_{U}\in \underline{K}^{U}}$ is compatible
with restrictions to faces. Let $\Delta$ and $\Delta'$ be two simplices of $K$
such that $s_{j}\prec \Delta'\prec \Delta$ for some $j\in \{1,\dots,e\}$. Then, $\Delta'_{U}\cap
V\subset \Delta_{U}\cap V$ and, since
$\gamma$ is piecewise smooth, we have
$$\gamma_{\Delta'_{U\cap V}}(x)=(\gamma_{\Delta_{U\cap V}})_{/\Delta'_{U\cap V}}(x)$$
for each $x\in \Delta'_{U\cap V}$. Hence, if $x\in \Delta'_{U\cap V}$,
$$\xi_{\Delta'_{U}}(x)=-\psi(x)\
\gamma_{\Delta'_{U\cap V}}(x)=-\psi(x)(\gamma_{\Delta_{U\cap V}})_{/\Delta'_{U\cap V}}(x)=
(\xi_{\Delta_{U}})_{/\Delta'_{U}}(x)$$ If $x\in
\Delta_{U}\cap(M\setminus \supp \psi)$ we have
that
$\xi_{\Delta_{U}}(x)=(\xi_{\Delta_{U}'})_{\Delta_{U}}(x)=0$.
Hence, the differential form $(\xi_{\Delta_{U}})_{\Delta_{U}\in \underline{K}^{U}}$ is a
piecewise smooth form on $\cA_{U}$. Analogously, we define a piecewise smooth
form $(\eta_{\Delta_{V}})_{\Delta_{V}\in \underline{K}^{V}}\in \Omega^{p}_{ps}(\cA_{V};V)$ by
$$\eta_{\Delta_{V}}(x)=\begin{cases}

-\varphi(x)\ \gamma_{\Delta_{U\cap V}}(x) & \text{if $x\in
\Delta_{V}\cap U$}
\\\\\
\ 0_{x}\in (\cA^{!!}_{\Delta_{V}})_{x} & \text{if $x\in
\Delta_{V}\cap(M\setminus \supp \varphi)$}
\end{cases}$$
and we have that, for each $x\in \Delta_{U\cap V}\in \underline{K}^{U\cap V}$,
$$(\eta_{\Delta_{V}})_{/\Delta_{U\cap V}}(x)-(\xi_{\Delta_{U}})_{/\Delta_{U\cap V}}(x)=\gamma_{\Delta_{U\cap V}}(x)$$ that is,
$$\pi\big((\xi_{\Delta_{U}})_{\Delta_{U}\in \underline{K}^{U}},(\eta_{\Delta_{V}})_{\Delta_{V}\in \underline{K}^{V}}\big)=(\gamma_{\Delta_{U\cap V}})_{\Delta_{U\cap V}\in \underline{K}^{U\cap V}}$$ Hence, the result is proved. {\small $\square$}

\vspace{3mm}

The Mayer-Vietoris sequence in the piecewise context is the long sequence of cohomology corresponding to the short sequence shown in the previous proposition. If we put this short sequence together with the short exact sequence presented by Kubarski in the third section of \cite{kuki-dual} (smooth case), we obtain a commutative diagram of short exact sequences, in which the vertical maps are the restriction maps $\Psi$ above. This is our next proposition.

\vspace{3mm}

\begin{prop} Keeping the same hypothesis and notation as above, the diagram
\[\minCDarrowwidth15pt\begin{CD}
\{0\} @>>> \Omega^{p}(\cA;M) @>\lambda>> \Omega^{p}(\cA_{U};U)\oplus \Omega^{p}(\cA_{V};V) @>\mu>> \Omega^{p}(\cA_{U\cap V};U\cap V) @>>> \{0\}\\
@. @VV\Psi V @VV\Psi V @VV\Psi V @.\\
\{0\} @>>> \Omega^{p}(\cA;K) @>\delta>> \Omega^{p}_{ps}(\cA_{U};U)\oplus \Omega^{p}_{ps}(\cA_{V};V) @>\pi>> \Omega^{p}_{ps}(\cA_{U\cap V};U\cap V) @>>> \{0\}\\
\end{CD}\]
is commutative.
\end{prop}

\vspace{3mm}

We finalize this section with a proposition concerning the K$\ddot{\textrm{u}}$nneth isomorphism in a particular case of piecewise setting. Let $M$ be a smooth manifold, smoothly triangulated by a simplicial complex $K$, and $\s$ a simplex of $K$. Let $U=\textbf{\textrm{St}}(\s)$. Assume that $\mathfrak{g}$ is a real Lie algebra and consider the trivial Lie algebroid $\cA=TU\oplus (M\times \frak g)$ on $U$, which is identified to the Lie algebroid $\cA=TU\times \frak g$ by a strong isomorphism of Lie algebroids over $U$. For each $p\geq 0$, one has $\Omega^{p}(\frak g)=\bigwedge^{p}\frak g$.

\vspace{3mm}

\begin{prop} Keeping the same hypothesis and notations as above,
$$\textrm{H}_{ps}(\cA;U)\simeq \textrm{H}_{ps}(U)\otimes \textrm{H}(\mathfrak{g})$$
\end{prop}

\proofwp. As done before, we denote by $\underline{K}^{U}$ the set of all submanifolds $\Delta_{U}=U\cap \Delta$ for each simplex $\Delta\in K$ such that $\s$ is a face of $\Delta$. Consider the Lie algebroids morphisms $$\gamma_{\Delta_{U}}:T\Delta_{U}\times \frak g \longrightarrow T\Delta_{U}$$ and $$\pi_{\Delta_{U}}:T\Delta_{U}\times \frak g \longrightarrow \frak g$$ given by the projections on the first and second factors respectively ($\pi$ is a non strong Lie algebroid morphism over a constant map). We divide the proof in three parts.

\vspace{2mm}

Part 1. We are going to check that
$$\Omega^{\ast}_{ps}(U)\otimes \Omega^{\ast}(\mathfrak{g})\simeq \Omega^{\ast}_{ps}(\cA;U)$$
Let $\xi=(\xi_{\Delta_{U}})_{\Delta_{U}\in \underline{K}^{U}}\in \Omega^{\ast}_{ps}(U)$ and $\eta\in \Omega^{\ast}(\frak g)$. If $\Delta'$ and $\Delta$ are two simplices of $K$ such that $s\prec \Delta'\prec \Delta$, denote by $$(\varphi_{\Delta,\Delta'}^{T\Delta_{U}\times \frak g})^{!!}:(T\Delta_{U}\times \frak g)^{!!}_{\Delta'_{U}}\longrightarrow T\Delta_{U}\times \frak g$$ and $$(\varphi_{\Delta,\Delta'}^{T\Delta_{U}})^{!!}:(T\Delta_{U})^{!!}_{\Delta'_{U}}\longrightarrow T\Delta_{U}$$ the canonical maps (see the paragraph after proposition 1.1) induced from the diagrams
\begin{equation*}
\begin{array}{ccccccccccccccc}
\  &  & T\Delta_{U}\times \frak g & & & & & & & & \  &  & T\Delta_{U}\\
& & \mapd{\gamma} & & & & & & & & & & \mapd{\pi_{T\Delta_{U}}} \\
\\ \Delta'_{U} & \mapr{\varphi_{\Delta,\Delta'}} & \Delta_{U} & & & & & & & & \Delta'_{U} & \mapr{\varphi_{\Delta,\Delta'}} & \Delta_{U}
\end{array}
\end{equation*}
It is obvious that
$$\gamma_{\Delta_{U}}\circ (\varphi_{\Delta,\Delta'}^{T\Delta_{U}\times \frak g})^{!!}=(\varphi_{\Delta,\Delta'}^{T\Delta_{U}})^{!!}\circ \gamma_{\Delta'_{U}}\ \ \ \textrm{and}\ \ \ \pi_{\Delta_{U}}\circ (\varphi_{\Delta,\Delta'}^{T\Delta_{U}\times \frak g})^{!!}=\pi_{\Delta'_{U}}$$ so the equalities
$$(\gamma^{\ast}_{\Delta_{U}}\xi_{\Delta_{U}})_{/\Delta_{U}'}=\gamma^{\ast}_{\Delta_{U}'}\xi_{\Delta'_{U}}\ \ \ \textrm{and}\ \ \ (\pi^{\ast}_{\Delta_{U}}\eta)_{/\Delta_{U}'}=\pi^{\ast}_{\Delta_{U}'}\eta$$ hold too. These equalities show that the differential form $$(\gamma^{\ast}_{\Delta_{U}}\xi_{\Delta_{U}}\wedge \pi^{\ast}_{\Delta_{U}}\eta)_{\Delta_{U}\in \underline{K}^{U}}$$ belongs to $\Omega^{\ast}_{ps}(\cA;U)$. Hence, we can consider a map $$k_{ps}:\Omega^{\ast}_{ps}(U)\otimes \Omega^{\ast}(\frak g)\longrightarrow \Omega^{\ast}_{ps}(\cA;U)$$ such that $$k_{ps}(\xi\otimes \eta)=(\gamma^{\ast}_{\Delta_{U}}\xi_{\Delta_{U}}\wedge \pi^{\ast}_{\Delta_{U}}\eta)_{\Delta_{U}\in \underline{K}^{U}}$$ where $\xi=(\xi_{\Delta_{U}})_{\Delta_{U}\in \underline{K}^{U}}\in \Omega^{\ast}_{ps}(U)$ and $\eta\in \Omega^{\ast}(\frak g)$. This map is well defined. Now, we shall see that the map $k_{ps}$ is an isomorphism of differential graded algebras. Obviously, the map $k_{ps}$ is a morphism of graded algebras. For each $\Delta\in K$ such that $s\prec \Delta$, let $$k_{\Delta_{U}}:\Omega^{\ast}(\Delta_{U})\otimes \Omega(\frak g)\longrightarrow \Omega^{\ast}(T\Delta_{U}\times \frak g;\Delta_{U})$$ be the K$\ddot{\textrm{u}}$nneth isomorphism described by Kubarski in the sixth section of \cite{kuki-dual}. We have that, $$(k_{ps}(\xi\otimes \eta))_{\Delta_{U}}=\gamma^{\ast}_{\Delta_{U}}\xi_{\Delta_{U}}\wedge \pi^{\ast}_{\Delta_{U}}\eta=k_{\Delta_{U}}(\xi_{\Delta_{U}}\otimes \eta)$$ Therefore, if $\omega=\sum \xi\otimes \eta\in \Omega^{\ast}_{ps}(U)\otimes \Omega^{\ast}(\frak g)$ and $k_{ps}(\omega)=0$, then $(k_{ps}(\omega))_{\Delta_{U}}=0$ and so $$0=\big(k_{ps}(\sum \xi\otimes \eta)\big)_{\Delta_{U}}=k_{\Delta_{U}}\big(\sum (\xi_{\Delta_{U}}\otimes \eta)\big)$$ Hence $\omega=\sum (\xi_{\Delta_{U}}\otimes \eta)=0$ and, with this, we have checked that $k$ is injective. Take now
$\lambda=(\lambda_{\Delta_{U}})_{\Delta_{U}\in \underline{K}^{U}}\in \Omega^{\ast}_{ps}(\cA;U)$. We want to find $\omega\in \Omega^{\ast}_{ps}(U)\otimes \Omega^{\ast}(\frak g)$ such that $k_{ps}(\omega)=\lambda$. Since $k_{\Delta_{U}}$ is surjective, we can consider smooth forms $\xi_{j_{\Delta_{U}}}\in \Omega^{\ast}(\Delta_{U})$ and $\eta\in \Omega^{\ast}(\frak g)$ such that $$k_{\Delta_{U}}\big(\sum_{j}(\xi_{j_{\Delta_{U}}}\otimes \eta)\big)=\lambda_{\Delta_{U}}$$ Take then the form $\omega_{\Delta_{U}}=\sum_{j}(\xi_{j_{\Delta_{U}}}\otimes \eta)$. If $\Delta'$ and $\Delta$ are simplices of $K$ with $s\prec \Delta'\prec \Delta$, we have the equalities
$$k_{\Delta'_{U}}\big(\sum_{j}(\xi_{j_{\Delta_{U}}})_{\Delta_{U}'}\otimes \eta\big)=\sum_{j}k_{\Delta'_{U}}((\xi_{j_{\Delta_{U}}})_{\Delta_{U}'}\otimes \eta\big)=\sum_{j}\big(\gamma^{\ast}_{\Delta_{U}'}(\xi_{j_{\Delta_{U}}})_{\Delta_{U}'}\wedge \pi^{\ast}_{\Delta'}\eta\big)=(\ast)$$ and $$k_{\Delta'_{U}}\big(\sum_{j}(\xi_{j_{\Delta'_{U}}}\otimes \eta)\big)=\lambda_{\Delta'_{U}}=(\lambda_{\Delta_{U}})_{/\Delta'_{U}}=\big(k_{\Delta_{U}}(\sum_{j}(\xi_{j_{\Delta_{U}}}\otimes \eta)\big)_{/\Delta_{U}'}=$$ $$=\big(\sum_{j}(\gamma^{\ast}_{\Delta_{U}}(\xi_{j_{\Delta_{U}}})\wedge \pi^{\ast}_{\Delta}\eta)\big)_{/\Delta_{U}'}=\sum_{j}(\gamma^{\ast}_{\Delta_{U}}(\xi_{j_{\Delta}})\wedge \pi^{\ast}_{\Delta}\eta)\big)_{/\Delta_{U}'}=$$ $$=\sum_{j}(\gamma^{\ast}_{\Delta_{U}}((\xi_{j_{\Delta_{U}}})_{/\Delta'_{U}}))\wedge \pi^{\ast}_{\Delta'}\eta))=\sum_{j}\big(\gamma^{\ast}_{\Delta_{U}'}(\xi_{j_{\Delta_{U}}})_{\Delta'_{U}}\wedge \pi^{\ast}_{\Delta'}\eta\big)=(\ast)$$ Hence, $$k_{\Delta'_{U}}\big(\sum_{j}(\xi_{j_{\Delta_{U}}})_{\Delta_{U}'}\otimes \eta\big)=k_{\Delta'_{U}}\big(\sum_{j}(\xi_{j_{\Delta'_{U}}}\otimes \eta)\big)$$ and, since $k_{\Delta'_{U}}$ is bijective, $\sum_{j}(\xi_{j_{\Delta_{U}}})_{\Delta_{U}'}\otimes \eta=\sum_{j}(\xi_{j_{\Delta'_{U}}}\otimes \eta)$. Therefore, we can conclude that $(\xi_{j_{\Delta_{U}}})_{/\Delta_{U}'}=\xi_{j_{\Delta'_{U}}}$. Then, the form $\omega=(\omega_{\Delta_{U}})_{\Delta_{U}\in \underline{K}^{U}}$, in which $\omega_{\Delta_{U}}=\sum_{j}(\xi_{j_{\Delta_{U}}}\otimes \eta)$, belongs to $\Omega^{\ast}_{ps}(U)\otimes \Omega^{\ast}(\frak g)$. Obviously, $k_{ps}(\omega)=\lambda$ and then it is checked that $k_{ps}$ is an isomorphism of graded algebras.

\vspace{2mm}

Part 2. Next, we are going to check that $k_{ps}$ commutes with differential. For each $\Delta\in K$ such that $\s$ is a face of $\Delta$, the differentials on the complexes $\Omega^{\ast}_{ps}(\cA;U)$, $\Omega^{\ast}_{ps}(U)$ and $\Omega^{\ast}(\Delta_{U})$ are denoted by $d^{\cA}_{ps}$, $d^{U}_{ps}$ and $d^{\Delta_{U}}$ respectively. Let $\xi=(\xi_{\Delta_{U}})_{\Delta_{U}\in \underline{K}^{U}}\in \Omega^{\ast}_{ps}(U)$ and $\eta\in \Omega^{\ast}(\frak g)$. We have
$$(d^{\cA}_{ps}\circ k_{ps})(\xi\otimes \eta)=d^{\cA}_{ps}((\gamma^{\ast}_{\Delta_{U}}\xi_{\Delta_{U}}\wedge \pi^{\ast}_{\Delta_{U}}\eta)_{\Delta_{U}\in \underline{K}^{U}})=$$ $$=d^{\cA}_{ps}\big((\gamma^{\ast}_{\Delta_{U}}\xi_{\Delta_{U}})_{\Delta_{U}\in \underline{K}^{U}}\big)\wedge \pi^{\ast}\eta+(-1)^{deg\xi}(\gamma^{\ast}_{\Delta_{U}}\xi_{\Delta_{U}})_{\Delta_{U}\in \underline{K}^{U}}\wedge d^{\cA}_{ps}(\pi^{\ast}\eta)=$$
$$=\big(\gamma^{\ast}_{\Delta_{U}}(d^{\Delta_{U}}(\xi_{\Delta_{U}}))\big)_{\Delta_{U}\in \underline{K}^{U}}\wedge \pi^{\ast}\eta+(-1)^{deg \omega}\gamma^{\ast}\xi\wedge \pi^{\ast}(d_{\frak g}\eta)=$$
$$=k_{ps}((d^{U}_{ps}\xi)\otimes \eta)+(-1)^{deg \xi}k_{ps}(\xi\otimes d_{\frak g}\eta)=k_{ps}\circ \delta(\xi\otimes \eta)$$ This proves that $k_{ps}$ is an isomorphism of differential graded algebras.

\vspace{2mm}

Part 3. The isomorphism $k_{ps}$ above induces an isomorphism in cohomology. By applying the K$\ddot{\textrm{u}}$nneth theorem, we obtain $$H^{\ast}_{ps}(\cA;U)\simeq H^{\ast}(\Omega^{\ast}_{ps}(U)\otimes \Omega^{\ast}(\mathfrak{g}))\simeq H^{\ast}_{ps}(U)\otimes H^{\ast}(\frak g)$$ and the result is proved. {\small $\square$}


\begin{center}

\section{Main theorem}

\end{center}

\vspace{6mm}

Whitney in \cite{wity-git} and Sullivan in \cite{suli-inf} have shown that cohomologies obtained by using the cell structure of a space are isomorphic to the singular cohomology of the polytope. Therefore, those piecewise cohomologies also are isomorphic to the Rham cohomology, if the space is a cell smooth manifold. Based both in their work and in the work developed by Mackenzie (see \cite{makz-gad} and \cite{makz-lga}) as well by Kubarski (see \cite{kuki-dual} and \cite{kuki-chw}), we have claimed in the previous section that the Lie algebroid cohomology and
the piecewise smooth cohomology of a transitive Lie algebroid over a triangulated compact manifold are isomorphic. We proceed to the proof of this assertion for all Lie algebroids under these hypothesis.

\vspace{3mm}

Let $M$ be a smooth manifold, smoothly triangulated by a finite simplicial complex $K$, and $\cA$ a transitive Lie algebroid on $M$. Consider the restriction map $$\Psi:\Omega^{\ast}(\cA;M)\longrightarrow \Omega^{\ast}(\cA;K)$$ $$\omega \longrightarrow (\omega_{/\Delta})_{\Delta\in
K}$$ This section is then devoted to prove the following theorem.

\vspace{3mm}

\begin{thm} The map $\Psi$ induces an isomorphism in cohomology.
\end{thm}

\vspace{3mm}

The proof of this theorem involves, beyond theorems already mentioned before, the Steenrod lemma applied to the commutative diagram shown in the proposition 3.2 of the previous section. In fact, we shall be able to apply the Steenrod lemma if we know that the map $\Psi$ is a quasi-isomorphism for trivial Lie algebroids over stars. Therefore, our first step is to show that the main theorem holds for these trivial Lie algebroids.

\vspace{3mm}

\begin{prop} Let $M$ be a smooth manifold, smoothly triangulated by a simplicial complex $K$, and $\s$ a simplex of $K$. Let $U=\textbf{\textrm{St}}(\s)$ and denote by $\underline{K}^{U}$ the set of all submanifolds $\Delta_{U}=U\cap \Delta$ such that $\Delta\in K$ and $\s$ is a face of $\Delta$. Assume that $\mathfrak{g}$ is a real Lie algebra and consider the trivial Lie algebroid $\cA=TU\times \frak g$. Then,
the restriction map
$$\Psi:\Omega^{\ast}(\cA;U)\longrightarrow \Omega^{\ast}_{ps}(\cA;U)$$
$$\omega \longrightarrow (\omega_{/\Delta_{U}})_{\Delta_{U}\in \underline{K}^{U}}$$ induces an isomorphism in cohomology.
\end{prop}

\proofwp. By the K$\ddot{\textrm{u}}$nneth theorem for trivial Lie algebroids stated by Kubarski in \cite{kuki-dual}, we have that
$$\textrm{H}(\cA;U)\simeq \textrm{H}_{\textrm{dR}}(U)\otimes \textrm{H}(\mathfrak{g})$$
Next, we shall see that $\Psi$ induces an isomorphism in cohomology. Take the diagram
$$
\xymatrix{
\Omega^{\ast}(U)\otimes \Omega^{\ast}(\mathfrak{g})\ar[d]_{k}\ar[r]^{\lambda} & \Omega^{\ast}_{ps}(U)\otimes \Omega^{\ast}(\mathfrak{g})\ar[d]^{k_{ps}} \\
 \Omega^{\ast}(\cA;U)\ar[r]^{\Psi} & \Omega^{\ast}_{ps}(\cA;U)
}
$$
where $k_{ps}:\Omega^{\ast}_{ps}(U)\otimes \Omega^{\ast}(\frak g)\longrightarrow \Omega^{\ast}_{ps}(\cA;U)$ is the isomorphism defined in the proof of the proposition 3.5 of the previous section, $k$ is the K$\ddot{\textrm{u}}$nneth isomorphism described by Kubarski in the sixth section of \cite{kuki-dual} and $\lambda=\Phi\otimes \id$, in which $\Phi$ is the restriction map obtained from the Rham-Sullivan theorem for cell manifolds (see Sullivan \cite{suli-inf} or Whitney \cite{wity-git}). Obviously, the diagram is commutative and, by the de Rham-Sullivan theorem (see Sullivan \cite{suli-inf}), the map $\Phi$ induces an isomorphism in cohomology. Therefore, in cohomology, we have the commutative diagram
$$
\xymatrix{
H^{\ast}_{dR}(U)\otimes H^{\ast}(\mathfrak{g})\ar[d]_{\simeq} \ar[r]^{H(\lambda)} & H^{\ast}_{ps}(U)\otimes H^{\ast}(\mathfrak{g})\ar[d]_{\simeq} \\
H^{\ast}(\cA;U)\ar[r]^{H(\Psi)} & H^{\ast}_{ps}(\cA;U)
}
$$
Hence, $H(\Psi)$ is an isomorphism. {\small $\square$}

\vspace{3mm}

Next, we want to show that $\Psi$ induces an isomorphism in cohomology not only for trivial Lie algebroids defined over regular open subsets but for any arbitrary transitive Lie algebroid defined over a regular open subset. For that, we begin with a basic result which is a direct consequence from the functor homology.

\vspace{3mm}

\begin{prop} Let $M$ be a smooth manifold, smoothly triangulated by a simplicial complex $K$, $\s$ a simplex of $K$ and $U=\St(\s)$. Let $\cA$ and $\cB$ be two transitive Lie algebroids on $M$ and suppose there is an isomorphism of Lie algebroids between them. Then, the cohomology spaces $H_{ps}(\cA;U)$ and $H_{ps}(\cB;U)$ are isomorphic.
\end{prop}

\vspace{3mm}

\begin{prop} Let $M$ be a smooth manifold, smoothly triangulated by a simplicial complex $K$, $\s$ a simplex of $K$ and $U=\St(\s)$. Denote by $\underline{K}^{U}$ the set of all submanifolds $\Delta_{U}=U\cap \Delta$ such that $\Delta\in K$ and $\s$ is a face of $\Delta$. Let $\cA$ be a transitive Lie algebroid on $U$.  Then, the morphism
\vspace{3mm}
$$\Psi:\Omega^{\ast}(\cA;U)\longrightarrow \Omega^{\ast}_{ps}(\cA;U)$$ $$\omega \longrightarrow (\omega_{/\Delta_{U}})_{\Delta_{U}\in \underline{K}^{U}}$$
induces an isomorphism in cohomology.
\end{prop}

\proofwp. Since $U$ is contractible, $\cA$ is isomorphic to the trivial Lie algebroid $\cB=TU\times \frak g$ on $U$, in which $\frak g$ is the fibre type
$\ker \gamma$. We conclude the result by the commutativity of the diagram
\begin{equation*}
\begin{array}{ccccc}
\Omega^{p}(\cA;U) & \mapr{} & \Omega^{p}_{ps}(\cA;U)\\
 \mapd{} & & \mapd{} \\
 \Omega^{p}(TU\times \frak g) & \mapr{\Psi} & \Omega^{p}_{ps}(TU\times \frak g)\\
\end{array}
\end{equation*}
and applying the proposition 4.3 to the vertical maps and the proposition 4.2 to the map $\Psi$ displayed on the diagram. {\small $\square$}

\vspace{3mm}

\textit{Proof of the theorem}. We shall prove the result by induction on the number of vertices of the simplicial complex $K$. Suppose then that $v_{0}$, $\dots$, $v_{n}$ is the family of all vertices of $K$. If $K$ has only one vertex, the result is trivial. Suppose we have established the result for all $l< n$. It is well known that the family $\{\St (v_{j}):j\in \{0,\dots, n\}\}$ is an open covering of $M$. Consider the open subsets $U=\bigcup^{n-1}_{j=0}\St (v_{j})$ and $V=\St (v_{n})$ of $M$. We have that
$$U\cap V=\big(\bigcup^{n-1}_{j=0}\St(v_{j})\big)\cap \St(v_{n})=\bigcup^{n-1}_{j=0}\big(\St(v_{j})\cap \St(v_{n})\big)=\bigcup_{j}\St ([v_{j},v_{n}])$$
where the union is taken over all indexes $j$ such that the vertices $v_{j}$ and $v_{n}$ generates a simplex of $K$
(otherwise the intersection $\St(v_{j})\cap \St(v_{n})$ is empty) and $\big[v_{j},v_{n}\big]$ denotes the closed simplex generated by the vertices $v_{j}$ and $v_{n}$. We have seen that
\[\minCDarrowwidth15pt\begin{CD}
\{0\} @>>> \Omega^{p}(\cA;M) @>\lambda>> \Omega^{p}(\cA_{U};U)\oplus \Omega^{p}(\cA_{V};V) @>\mu>> \Omega^{p}(\cA_{U\cap V};U\cap V) @>>> \{0\}\\
@. @VV\Psi V @VV\Psi V @VV\Psi V @.\\
\{0\} @>>> \Omega^{p}(\cA;K) @>\delta>> \Omega^{p}_{ps}(\cA_{U};U)\oplus \Omega^{p}_{ps}(\cA_{V};V) @>\pi>> \Omega^{p}_{ps}(\cA_{U\cap V};U\cap V) @>>> \{0\}\\
\end{CD}\] is a commutative diagram of short exact sequences. The map $\Psi$ on the right side is quasi-isomorphism by induction. The map $\Psi$ on the middle is quasi-isomorphism by induction and the previous proposition. By the Steenrod lemma, the map $\Psi$ on the left side is also a quasi-isomorphism.
{\small $\square$}

\vspace{3mm}

From the main theorem, we easily infer that the piecewise smooth cohomology of a combinatorial compact manifold does not depend on the triangulation used, that is, for any simplicial division of the simplicial complex, the piecewise smooth cohomology spaces of both combinatorial manifolds remains isomorphic. Precisely, this statement is our next proposition.

\vspace{3mm}

\begin{cor} Let $M$ be a smooth manifold, smoothly triangulated by a simplicial complex $K$, and $\cA$ a transitive Lie algebroid on $M$. Let $L$ be other simplicial complex and assume that $L$ is a subdivision of $K$. Then, the piecewise smooth cohomology $H(\cA;K)$ is isomorphic to the piecewise smooth cohomology $H(\cA;L)$. Thus, the morphism from $\Omega^{\ast}(\cA;K)$ to $\Omega^{\ast}(A;L)$ which induces that isomorphism in cohomology is also given by restriction of forms.
\end{cor}

\proofwp. The result follows from the commutativity of the following
diagram

$$
\xymatrix{
& \Omega^{\ast}(\cA;M) \ar[dl]_{\Psi} \ar[dr]^{\Psi} & & \\
\Omega^{\ast}(\cA;K) \ar[rr]^{\Phi} & & \Omega^{\ast}(\cA;L)
}
$$
where $\Phi$ is also given by restriction. {\small $\square$}

\vspace{3mm}

\Addresses

\end{document}